%% file: Bryan-Kool.tex
\documentclass[12pt]{amsart}

\title[DT invariants of local elliptic surfaces]{Donaldson-Thomas
invariants of local elliptic surfaces via the topological vertex}
\author{Jim Bryan and Martijn Kool} \date{\today} \address{
Department of Mathematics\\
University of British Columbia \\
Room 121, 1984 Mathematics Road  \\
Vancouver, B.C., Canada V6T 1Z2  
}

\address{
Mathematical Institute \\
Utrecht University \\
Room 502, Budapestlaan 6  \\
3584 CD Utrecht, The Netherlands  
}

\usepackage{diagrams}
\usepackage{mathtools}

\usepackage{tikz}
\usepackage{pgfplots}
\usepackage{tikz}
\usepackage{verbatim}

\usetikzlibrary{calc,3d}

\definecolor{tealgreen}{HTML}{1B9E77}
\definecolor{orange}{HTML}{D95F02}
\definecolor{purple}{HTML}{7570B3}
\definecolor{pink}{HTML}{E7298A}
\definecolor{grassgreen}{HTML}{66A61E}
\definecolor{goldyellow}{HTML}{E6AB02}
\definecolor{brown}{HTML}{A6761D}
\definecolor{devilgray}{HTML}{666666}

\usepackage{amsmath}
\usepackage{verbatim}
\usepackage{amsmath,amsthm,amsfonts}
\usepackage{amssymb}
\usepackage{times}
\usepackage[colorinlistoftodos]{todonotes}

\newtheorem{theorem}{Theorem}
\newtheorem{proposition}[theorem]{Proposition}
\newtheorem{conjecture}[theorem]{Conjecture}
\newtheorem{lemma}[theorem]{Lemma}
\newtheorem{corollary}[theorem]{Corollary}
\newtheorem{convention}{Convention}[theorem]

\theoremstyle{definition}

\newtheorem{def-theorem}[theorem]{Definition-Theorem}
\newtheorem{remark}[theorem]{Remark}
\newtheorem{definition}[theorem]{Definition}

\newcommand{\CC} {\mathbb{C}}          
\newcommand{\NN} {\mathbb{N}}		
\newcommand{\RR} {\mathbb{R}}		
\newcommand{\ZZ} {\mathbb{Z}}		
\newcommand{\PP} {\mathbb{P}}

\renewcommand{\O}{\mathcal{O}}

\newcommand{\sfV}{\mathsf{V}}
\newcommand{\sfVtilde}{\widetilde{\mathsf{V}}}
\newcommand{\Pic}{\mathrm{Pic}}

\newcommand\Into{\ensuremath{\lhook\joinrel\relbar\joinrel\rightarrow}}
\newcommand\INTO{\ar@{^{(}->}[r]}

\newcommand{\Hom}{\operatorname{Hom}}

\newcommand{\Ext}{\operatorname{Ext}}
\newcommand{\Ker}{\operatorname{Ker}}

\newcommand{\id}{\operatorname{id}}

\newcommand{\Sym}{\operatorname{Sym}}

\newcommand{\Hilb}{\operatorname{Hilb}}
\newcommand{\Quot}{\operatorname{Quot}}
\newcommand{\Ann}{\operatorname{Ann}}
\newcommand{\Tot}{\operatorname{Tot}}
\newcommand{\DT}{\mathsf{DT}}
\newcommand{\CM}{\operatorname{CM}}
\newcommand{\Var}{\operatorname{Var}}

\newcommand{\Supp}{\operatorname{Supp}}
\newcommand{\Spec}{\operatorname{Spec}}
\newcommand{\sm}{\operatorname{sm}}
\newcommand{\sing}{\operatorname{sing}}
\newcommand{\conn}{\operatorname{conn}}
\newcommand{\F}{\mathcal{F}}

\newcommand{\boldx}{\boldsymbol{x}}
\newcommand{\boldy}{\boldsymbol{y}}
\newcommand{\bolda}{\boldsymbol{a}}
\newcommand{\boldb}{\boldsymbol{b}}
\newcommand{\boldlambda}{\boldsymbol{\lambda }}
\newcommand{\boldmu}{\boldsymbol{\mu }}

\renewcommand{\emptyset}{\varnothing}
\renewcommand{\hat}{\widehat}
\newcommand{\DThat}{\widehat{\DT}}

\newcommand{\fiber}{\mathsf{fib}}
\newcommand{\PCP}{\mathsf{PCP}}
\newcommand{\PFP}{\mathsf{PFP}}
\newcommand{\half}{\frac{1}{2}}
\newcommand{\red}{\mathrm{red}}
\newcommand{\bx}{\square}

\newcommand{\mujprime}{\mu^{(j)\prime}}
\newcommand{\length}{\operatorname{leng}}
\newcommand{\presectionspace}{\vspace{0.2cm}} 

\newcommand{\SubSecSpace}{$\,$\vspace{0.2cm}\par } 

\begin{document}

\begin{abstract}

We compute the Donaldson-Thomas invariants of a local elliptic surface
with section. We introduce a new computational technique which is a
mixture of motivic and toric methods. This allows us to write the
partition function for the invariants in terms of the topological
vertex. Utilizing identities for the topological vertex proved in
\cite{Bryan-Kool-Young}, we derive product formulas for the partition
functions. The connected version of the partition function is written
in terms of Jacobi forms.  In the special case where the elliptic
surface is a K3 surface, we get a derivation of the
Katz-Klemm-Vafa formula for primitive curve classes which is
independent of the computation of Kawai-Yoshioka.
\end{abstract}

\maketitle 



\presectionspace
\section{Introduction}

Let $p : S \rightarrow B$ be a non-trivial elliptic surface over a
complex smooth projective curve $B$. We assume $p$ has a section and
all singular fibers are irreducible rational nodal curves. 

In this paper, we study the Donaldson-Thomas (DT) invariants of $X =
\Tot(K_S)$, i.e.~the total space of the canonical bundle $K_S$. This
is a non-compact Calabi-Yau threefold. Let $\beta$ be an effective
curve class on $S$. Consider the Hilbert scheme
$$
\Hilb^{\beta,n}(X) = \{ Z \subset X \ : \ [Z] = \beta, \ \chi(\O_Z) = n\}
$$
of proper subschemes $Z \subset X$ with homology class $\beta$ and holomorphic Euler characteristics $n$. The DT invariants of $X$ can be defined as
$$
\DT_{\beta,n}(X) := e(\Hilb^{\beta,n}(X), \nu) := \sum_{k \in \ZZ} k \ e(\nu^{-1}(k)),
$$
where $e(\cdot)$ denotes topological Euler characteristic and $\nu :
\Hilb^{\beta,n}(X) \rightarrow \ZZ$ is Behrend's constructible
function \cite{Behrend-micro}. We also consider an unweighted Euler
characteristic version of these invariants
$$
\DThat _{\beta,n}(X) := e(\Hilb^{\beta,n}(X)).
$$
We choose a section $B \subset S$ and focus on the primitive classes
$\beta = B + dF$, where $B$ is the class of the chosen section and $F$
the class of the fiber. We define the partition functions by
\begin{align*}
\DThat (X)& = \sum_{d = 0}^{\infty } \sum_{n \in \ZZ} \DThat _{B+dF,n}(X) p^n q^d,\\
\DT (X)& = \sum_{d =0}^{\infty } \sum_{n \in \ZZ} \DT_{B+dF,n}(X) y^n q^d.
\end{align*}

We also consider the partition functions for the invariants for
multiples of the fiber class
\begin{align*}
\DThat_{\fiber}(X)& = \sum_{d= 0}^{\infty} \sum_{n \in \ZZ} \DThat_{dF,n}(X) p^{n} q^{d} ,\\
\DT_{\fiber}(X)& = \sum_{d= 0}^{\infty} \sum_{n \in \ZZ} \DT_{dF,n}(X) y^{n} q^{d} .
\end{align*}

The main results of this paper are closed product formulas for the
partition functions $\DThat (X)$ and $\DThat_{\fiber}(X)$. Assuming a
general conjecture about the Behrend function, we also determine $\DT
(X)$ and $\DT_{\fiber}(X)$.

We use the notation
\[
M(p,q) = \prod_{m=1}^{\infty} (1-p^{m}q)^{-m}
\]
and the shorthand $M(p)=M(p,1)$.

\begin{theorem}\label{thm: main thm -- formulas for DT and DTfiber}
Let $e(S)$ and $e(B)$ denote the topological Euler characteristic of
the elliptic surface and the base. Then
\begin{align*}
\DThat (X) &= \left \{M(p)\prod_{d=1}^{\infty}
\frac{M(p,q^{d})}{(1-q^{d})} \right\}^{e(S)}
\left\{\frac{1}{(p^{\half}-p^{-\half})}\prod_{d=1}^{\infty}\frac{(1-q^{d})}{(1-pq^{d})(1-p^{-1}q^{d})}
\right\}^{e(B)} \\
\DThat_{\fiber}(X) &= \left\{M(p)\prod_{d=1}^{\infty}M(p,q^{d})
\right\}^{e(S)} \left\{\prod_{d=1}^{\infty}\frac{1}{(1-q^{d})}
\right\}^{e(B)}.
\end{align*}
\end{theorem}

The formula for $\DThat_{\fiber}(X)$ was previously proved using
wall-crossing methods by Toda\footnote{After applying the PT/DT
correspondence \cite{Bridgeland-PTDT}, this is essentially 
\cite[Thm~6.9]{Toda-2012-Kyoto}.}.

The ratio $\DThat (X)/\DThat_{\fiber}(X)$ can be considered as the
generating function for the connected invariants in the classes
$B+dF$. This series has a particularly nice form and can be written in
terms of classical Jacobi forms.  Consider the Dedekind eta function
and the Jacobi theta function
\begin{align*}
\eta &= q^{\frac{1}{24}} \prod_{k=1}^{\infty}(1-q^k), \\
\Theta &= 
(p^{\frac{1}{2}} - p^{-\frac{1}{2}}) \prod_{k=1}^{\infty} \frac{(1-p q^k) (1-p^{-1} q^k)}{(1-q^k)^2}.
\end{align*}

\begin{corollary} The partition function of the connected invariants is given as follows
\[
\frac{\DThat (X)}{\DThat_{\fiber}(X)}=
\left(q^{-\frac{1}{24}}\eta  \right)^{-e(S)}\Theta^{-e(B)}.
\]
\end{corollary}

In the case where $S \rightarrow \PP^1$ is an elliptically fibered
K3 surface, the above series specializes (up to a factor of $q$) to
the reciprocal of $\eta^{24}\Theta^{2}$, the unique Jacobi cusp form
of weight 10 and index 1. This is the Jacobi form appearing in the
well-known Katz-Klemm-Vafa formula. In order to obtain the KKV 
formula, we require the connected series, because $X$ is non-compact. 

Our result provides a new derivation of the KKV formula for primitive
classes\footnote{At least for the Euler characteristic version of the
Donaldson-Thomas invariants. For the Behrend function weighted
Donaldson-Thomas invariants, we require Conjecture~\ref{conj: Behrend
fnc conj}, see Theorem~\ref{thm: DT(X) assuming the Behrend function
conjecture} below.}. The KKV formula was proved in all curve classes
in \cite{Pandharipande-Thomas-KKV}. The appearance of the Jacobi form
$\eta^{24} \Theta^{2}$ in previous proofs of the KKV formula 
\cite{Maulik-Pandharipande-Thomas, Pandharipande-Thomas-KKV}
ultimately comes from the calculation of Euler characteristics 
of relative Hilbert schemes of points on curves on K3 by Kawai-Yoshioka 
\cite{Kawai-Yoshioka}. Our derivation of the KKV formula is the first that 
does not depend on the Kawai-Yoshioka formula.

Our results can be extended to apply to the usual (Behrend function
weighted) Donaldson-Thomas invariants if we assume a general
conjecture that we formulate in Section~\ref{sec: Behrend}. Our
conjecture relates the Behrend function at subschemes with embedded
points to the value of the Behrend function at the underlying
Cohen-Macaulay subscheme and may be of independent interest.

\begin{theorem}\label{thm: DT(X) assuming the Behrend function conjecture}
Assume that Conjecture~\ref{conj: Behrend fnc conj} holds, then
\[
\DT (X) = (-1)^{\chi (\O_{S})} \DThat (X)
\]
and
\[
\DT_{\fiber } (X) =  \DThat_{\fiber } (X)
\]
under the change of variables
\[
y=-p.
\]
\end{theorem}

A similar phenomenon to the above is known to hold when $X$ is a toric
Calabi-Yau threefold.

The method of computation that we introduce in this paper has
been applied to other elliptically fibered geometries. Indeed, it has
found applications to the calculation of DT generating functions on
$\textrm{K3} \times E$, where $E$ is an elliptic curve \cite{Bryan-K3xE} and
abelian threefolds \cite{BOPY}, and is expected to apply to $(\textrm{K3}\times
E)/G$ where $G$ is a finite group acting symplectically on each
factor. 

Although the geometry under consideration is not toric, we
combine $\CC^*$-localization, motivic methods, and
$(\CC^{*})^{3}$-localization to end up with expressions that only
depend on the topological vertex $\sfV_{\lambda\mu\nu}$, and the
topological Euler characteristics $e(B)$, $e(S)$. The outline of
our method is as follows:
\begin{itemize}
\item The $\CC^{*}$-action on $X$ induces an action on $\Hilb (X)$
whose Euler characteristic localizes to the $\CC^{*}$-fixed locus.  In
Section~\ref{sec: reduction to thickened comb curves} we show that any
$\CC^{*}$-invariant subscheme has a maximal Cohen-Macaulay subscheme
which is a curve of a special form which we call a partition thickened
comb curve (Definition~\ref{defn: partition thickened comb
curve}). This curve is determined by data consisting of points
$x_{i}\in B$ labelled by integer partitions $\lambda^{(i)}$. This
gives rise to a constructible morphism $\rho$ to $\Sym B$ taking the
value $\sum_{i} |\lambda^{(i)}|x_{i}$ on such a curve (see
Theorem~\ref{thm: C* invariant curves are partition thickened comb
curves with embedded points}).

\item In Section~\ref{sec: pushforward to sym prod}, we push forward
the Euler characteristic measure to $\Sym B$ via the map $\rho$. We
show that $\rho_{*}(1)$, the push-forward measure, has nice
multiplicative properties that allow us to compute the weighted
Euler characteristic over $\Sym B$ using a general result about
symmetric products (Lemma~\ref{lem: formula for euler char of sym
products}). 
\item To compute the push-forward measure $\rho_{*}(1)$ explicitly, we
must compute the Euler characteristics of the fibers of $\rho$. These
fibers are strata in the Hilbert scheme parameterizing subschemes
whose maximal Cohen-Macaulay subscheme is a fixed partition thickened
comb curve $C$. 
Here it is useful to switch from Hilbert schemes with fixed $C$ 
to Quot schemes of its ideal sheaf $I_C$. We 
introduce a further stratification of these Quot schemes by specifying the 
set theoretic support of the quotients. This allows us to write these Quot 
schemes as products (in $K$-theory) of Quot schemes 
of $I_C$ where the quotient is supported \emph{only} on one of the nodes of $C_{\red}$, 
or only on one of the components of $C_{\red}$ (minus the nodes), or only on the 
complement of $C$. The Quot scheme of quotients supported at a node of 
$C_{\red}$ can be expressed (in $K$-theory) as a Quot scheme of a partition 
thickened comb curve (determined by $C$) on $\CC^3$. Similarly, after further 
push-forwards to further symmetric products, we express the Euler characteristics 
of all other Quot schemes in terms of Euler characteristics of Quot 
schemes on $\CC^3$ as well (see Section~\ref{sec: restriction to formal nghds}). 
\item The Quot schemes of $\CC^3$ of the previous step all carry a natural $T=\left(\CC^{*}
\right)^{3}$-action. $T$-localization then allows us to write their
Euler characteristics in terms of the topological vertex (see
Section~\ref{sec: reduction to the vertex}). 
\item Finally, using the trace formulas for the topological vertex
proved in \cite{Bryan-Kool-Young}, we write our expression for
$\DThat (X)$ as the closed product formula given in Theorem~\ref{thm:
main thm -- formulas for DT and DTfiber}. 
\end{itemize}

Our proof of Theorem~\ref{thm: DT(X) assuming the Behrend
function conjecture} requires Theorem~\ref{thm: Ext computation}, an
involved computation of $\Ext^{1}_{0}(I_{C},I_{C})$ for partition
thickened comb curves $C$. The proof of Theorem~\ref{thm: Ext
computation} occupies most of Section~\ref{sec: smoothness and
deformations} and while technical in nature, the method we introduce
(again a mixture of formally local toric methods and global geometry)
may be of independent interest to the experts.

\subsection{Acknowledgments}
\SubSecSpace

We would like to thank Tom Graber, Paul Johnson, Manfred
Lehn, Oliver Leigh, Davesh Maulik, Georg Oberdieck, Rahul
Pandharipande, J{\o}rgen Rennemo, Bal\'azs Szendr\H{o}i, Richard
Thomas, Qizheng Yin, and Benjamin Young for helpful conversations.
We particularly thank the anonymous referee, whose comments led us to 
replace our previous approach (involving formal schemes and fpqc covers) 
by a much clearer stratification technique, which simply keeps track of the 
location of the embedded points with respect to $C_\red$ (Section \ref{sec: restriction to formal nghds}). 
This improved the paper significantly.

Part of this work was done during the semester \emph{Enumerative geometry 
of moduli spaces of sheaves in low dimension} at CIB/EPFL, who 
provided support and excellent working conditions.



\presectionspace
\section{Definitions, notation, and conventions} \label{defnotcon}

Let $p : S \rightarrow B$ be an elliptic surface over a smooth
projective curve $B$. We assume:
\begin{enumerate}
\item $S$ is a non-trivial fibration,
\item $p$ has a section $B \subset S$,
\item all singular fibers of $p$ are irreducible rational nodal curves. 
\end{enumerate}
We note that the number of singular fibers is equal to $e(S)$.

We write $F_x$ for the fiber $p^{-1}(x)$ over a  point $x \in
B$. We choose a section $B \subset S$ and denote its class in $H_2(S)$
by $B$ as well. We denote the class of the fiber by $F \in H_2(S)$.

Let $X = \Tot(K_S)$ be the total space of the canonical bundle $K_S$. 
For brevity, we define
\begin{align*}
\Hilb^{d,n}(X) &:=\Hilb^{B+dF,n}(X), \\
\DThat _{d,n}(X) &:= \DThat_{B+dF,n}(X).
\end{align*}

Since we are dealing with generating functions and our calculations
involve motivic methods on the Hilbert schemes, it is useful to
introduce the following notation. We define
\begin{align*}
\Hilb^{d,\bullet}(X) := \sum_{n \in \ZZ} \Hilb^{d,n}(X) \, p^n,
\end{align*}
where we view the right hand side as a formal Laurent series whose
coefficients are elements in the Grothendieck ring of varieties,
i.e. $K_0(\Var_{\CC})(\!(p)\!)$.

\begin{convention}\label{conv: bullet convention}
When an index is replaced by a bullet, we will multiply by the
appropriate variable and sum over the index. We regard the result as a
formal power (or Laurent) series whose coefficients lie in
$K_{0}(\Var_{\CC})$ and we extend operations of the Grothendieck group
(addition, multiplication, Euler characteristic) to the series in the
obvious way.
\end{convention}

For example
\[
\Hilb^{\bullet ,\bullet}(X) = \sum_{d=0}^{\infty}\sum_{n\in \ZZ}
\Hilb^{d,n}(X) q^{d}p^{n}\in K_{0}(\Var_{\CC})(\!(p)\!)[[q]],
\]
so that we can write
\[
\DThat (X) = e(\Hilb^{\bullet ,\bullet}(X)).
\]

It is notationally convenient to treat an Euler characteristic
weighted by a constructible function as a Lebesgue integral, where the
measurable sets are constructible sets, the measurable functions are
constructible functions, and the measure of a set is given by its
Euler characteristic. In this language we have
\[
\DThat_{d ,n}(X) =  \int_{\Hilb^{d ,n}(X)} 1\, de, \quad \quad
\DT_{d ,n}(X) = \int_{\Hilb^{d ,n}(X)} \nu \, de
\]
and following the bullet convention we have
\[
\DThat(X) =  \int_{\Hilb^{\bullet ,\bullet }(X)} 1\, de, \quad \quad
\DT(X) = \int_{\Hilb^{\bullet ,\bullet }(X)} \nu \, de.
\]

We will also need notation for subsets of the Hilbert scheme which
parameterize those subschemes obtained by adding embedded points
and/or zero dimensional components to
some fixed Cohen-Macaulay curve.

\begin{definition}\label{defn: Hilb(U,C)}
Let $C\subset X$ be a (not necessarily reduced)
Cohen-Macaulay subscheme of dimension 1. Consider the Hilbert scheme 
of subschemes $Z \subset X$ of class $[Z] = [C] \in H_2(X)$ and 
$\chi(\O_Z) = \chi(\O_C)+n$. Inside this Hilbert schemes, we define the 
following closed subset
\begin{align*}
\Hilb^{n}(X,C) = \{Z\subset X \text{ such that }C\subset Z\text{ and
$I_{C}/I_{Z}$ has finite length $n$} \}.
\end{align*}
\end{definition}

Once the Cohen-Macaulay curve $C \subset X$ is fixed, it is useful to work with the Quot scheme $\Quot_X^n(I_C)$ of zero dimensional quotients of $I_C$ of length $n$. We have the following lemma.
\begin{lemma} \label{HilbtoQuot}
The following equality holds in $K_0(\Var_{\CC})(\!(p)\!)$
\[
\Hilb^\bullet(X,C) = \Quot^\bullet_X(I_C).
\]
\end{lemma}
\begin{proof}
The universal quotient $I_{C \times \Quot^n_X(I_C)} \twoheadrightarrow Q$ has flat kernel $I_{\mathcal{Z}}$. This provides a flat family $\mathcal{Z} \subset X \times  \Quot^n_X(I_C)$ which gives a morphism to the Hilbert scheme. The kernel of a quotient $I_C \twoheadrightarrow Q$, where $Q$ is zero dimensional of length $n$, is an ideal sheaf $I_Z \subset I_C$ satisfying
$$
n = \chi(Q) = \chi(I_C / I_Z) = \chi(\O_Z) - \chi(\O_C).
$$
Every $\CC$-valued point of $\Hilb^n(X,C)$ arises from a quotient $I_Z \twoheadrightarrow Q$ in this way. This gives a geometric bijection $\Quot^n_X(I_C) \rightarrow \Hilb^n(X,C)$ from which the lemma follows.
\end{proof}

\presectionspace
\section{Reduction to partition thickened comb curves}\label{sec: reduction to thickened comb curves}

The action of $\CC^*$ on the fibers of $X$ lifts to the moduli space
$\Hilb^{d,\bullet}(X)$. Therefore
$$
\int_{\Hilb^{d,\bullet}(X)} 1 \, de = \int_{\Hilb^{d,\bullet}(X)^{\CC^*}} 1 \, de.
$$

The main result of this section is a classification of the subschemes
parameterized by $\Hilb^{d,n}(X)^{\CC^{*}}$, namely the
$\CC^{*}$-invariant subschemes.  We find that the maximal
Cohen-Macaulay subscheme of a $\CC^{*}$-invariant subscheme is
determined by a point in $\Sym^{d}(B)$ along with some discrete data
(a collection of integer partitions). We begin with some notation.

\begin{definition}\label{defn: comb curves} 
Let $T=\Tot(K_{S}|_{B})$ and let $p:X\to T$ be the elliptic fibration
induced by the elliptic fibration $p:S\to B$. We say that a subscheme $C\subset
X$ is a \textbf{comb curve} if $C=B\cup p^{-1}(Z) $ where $Z\subset T$
is a zero dimensional subscheme which is set-theoretically supported on
$B$.
\end{definition}

Let $\lambda =(\lambda_{1}\geq \dotsb \geq \lambda_{l}) $ be an
integer partition. Then $\lambda$ determines a zero dimensional subscheme
$Z_{\lambda}\subset \Spec \CC [[r,s]]$ given by the monomial ideal
\begin{equation}\label{eqn: monomial ideal I given by a partition lambda}
I_{\lambda}=(r^{\lambda_{1}},r^{\lambda_{2}}s,\dotsc
,r^{\lambda_{l}}s^{l-1},s^{l}).
\end{equation}
In terms of $\lambda$ as a Young diagram, we note $(\rho ,\sigma)\in
\lambda$ if and only if $r^{\rho}s^{\sigma}\notin I_{\lambda}.$

\begin{definition}\label{defn: partition thickened comb curve}
Let $C=B\cup p^{-1}(Z)$ be a comb curve, let $x_{1},\dotsc ,x_{n}\in
B\subset T$ be the points where $Z$ is supported, an let
$(r_{i},s_{i})$ be formal local coordinates on $T$ about each point
$x_{i}$ so that $s_{i}$ vanishes on $S\cap T$ and $r_{i}$
vanishes on $R_{i}\cap T$ where $R_{i}=\Tot (K_{S}|_{F_{x_{i}}})$. We
say that $C$ is a \textbf{partition thickened comb curve} if there
exists partitions $\lambda^{(1)},\dotsc ,\lambda^{(n)}$ such that $Z$
is given by $Z_{\lambda^{(i)}}$ in the local coordinates
$(r_{i},s_{i})$ about $x_{i}$.
We denote such a curve by
$B\cup_{i}\left(\lambda^{(i)}F_{x_{i}} \right)$.\footnote{Specifically, writing $r:=r_i$, $s:=s_i$, $\lambda = (\lambda_1 \geq \cdots \geq \lambda_l):=\lambda^{(i)}$, $x:=x_i$, $F:=F_{x_i}$, and $t$ for a coordinate at $x$ vanishing on $T$, the ideal of $B \cup \lambda F$ in a formal neighbourhood of $x$ is given by $(s,t) \cdot (r^{\lambda_{1}},r^{\lambda_{2}}s,\dotsc ,r^{\lambda_{l}}s^{l-1},s^{l}) \subset \CC[[r,s,t]]$.} We say that a
subscheme $Z\subset X$ is a \textbf{partition thickened comb curve
with points (PCP)} if the maximal Cohen-Macaulay subscheme
$Z_{\CM}\subset Z$ is a partition thickened comb curve, in other words,
$Z$ is obtained from a partition thickened comb curve by adding
embedded points and/or zero dimensional components. We denote by
\[
\Hilb^{d,n}_{\PCP}(X)\subset \Hilb^{d,n}(X)
\]
the locus in the Hilbert scheme parametrizing partition thickened
comb curves with points.
\end{definition}

In the next section it will be important to notationally distinguish
between singular and smooth fibers. See Figure~\ref{fig: drawing of
Cohen-Macaulay C^* fixed subscheme} for an illustration of a partition
thickened comb curve with smooth fibers $\{F_{x_{i}} \}$ thickened by
partitions $\{\lambda^{(i)} \}$ and nodal fibers  $\{F_{y_{j}} \}$ thickened by
partitions $\{\mu^{(j)} \}$.

\begin{figure}
\begin{tikzpicture}[
                    z  = {-15},
		    scale = 1]

\begin{scope}[yslant=-0.35,xslant=0]


\begin{scope} [canvas is yz plane at x=0]
\draw [black](0,0) rectangle (3,5);
\end{scope}
\begin{scope} [canvas is xz plane at y=0]
\draw [black](0,0) rectangle (4,5);
\end{scope}
\begin{scope} [canvas is xz plane at y=0.5]
\draw [black](0,0) rectangle (4,5);
\end{scope}
\draw [black](0,0) rectangle (4,3);

\begin{scope} [canvas is xz plane at y=3]
\draw [thick,blue](1,5) rectangle (1.2,4.5);
\draw [thick,blue](1.2,5) rectangle (1.4,4.5);
\draw [thick,blue](1,4.5) rectangle (1.2,4.0);
\draw [thick,blue](3,5) rectangle (3.2,4.5);
\draw [thick,blue](3.2,5) rectangle (3.4,4.5);
\draw [thick,blue](3,4.5) rectangle (3.2,4.0);
\draw [thick,blue](3.2,4.5) rectangle (3.4,4.0);
\draw [thick,blue](3,4.0) rectangle (3.2,3.5);
\end{scope}

\begin{scope} [canvas is xz plane at y=0.5]
\draw [thick,blue](1,5) rectangle (1.2,4.5);
\draw [thick,blue](1.2,5) rectangle (1.4,4.5);
\draw [thick,blue](1,4.5) rectangle (1.2,4.0);
\draw [thick,blue](3,5) rectangle (3.2,4.5);
\draw [thick,blue](3.2,5) rectangle (3.4,4.5);
\draw [thick,blue](3,4.5) rectangle (3.2,4.0);
\draw [thick,blue](3.2,4.5) rectangle (3.4,4.0);
\draw [thick,blue](3,4.0) rectangle (3.2,3.5);
\end{scope}

\draw [ultra thick,orange] 
                   (1  ,0   ,5)
to [out=90,in=-90] (1  ,0.6 ,5)
to [out=90,in=-90] (0.5,1.5 ,5)
to [out=90,in=-90] (1  ,2.4 ,5)
to [out=90,in=-90] (1  ,3   ,5);


\draw [ultra thick,orange] (0,0.5,5)--(4,0.5,5);

\draw [ultra thick,orange] 
                    (3   ,0   ,5) 
to [out=90,in=0]    (2.3 ,1.8 ,5) 
to [out=180,in=90]  (2   ,1.5 ,5) 
to [out=270,in=180] (2.3 ,1.2 ,5) 
to [out=0,in=270]   (3   ,3   ,5);

\node [above] at (3,3,3.5) {$\mu^{(j)} $};
\node [below] at (3,0,5) {$F_{y_{j}}$};
\node [right] at (3,0.2,5) {$y_{j}$};
\node [right] at (2.8,1.5,5) {$z_{j}$};
\node [above] at (1,3,4) {$\lambda^{(i)}$};
\node [below] at (1,0,5) {$F_{x_{i}}$};
\node [right] at (1,0.2,5) {$x_{i}$};
\node [left] at (0,2.0,5) {$S$};
\node [right] at (5,2.9,1.5) {$X=\Tot(K_{S})$};
\node [right] at (5,1.4,1.5) {$T=\Tot(K_{S}|_{B})$};
\node [right] at (5.25,2.3,1.5) {$p$};
\node [left] at (0,0.5,5) {$B$};

\draw [->] (5.25,2.75,1.5)--(5.25,1.75,1.5);

\begin{scope} [canvas is yz plane at x=4]
\draw [black](0,0) rectangle (3,5);
\end{scope}
\begin{scope} [canvas is xz plane at y=3]
\draw [black](0,0) rectangle (4,5);
\end{scope}
\draw [black](0,0,5) rectangle (4,3,5);
\draw [black,fill, opacity=0.1](0,0,5) rectangle (4,3,5);

\end{scope}
\end{tikzpicture}
\caption{A partition thickened comb curve
\[
C=B\cup_{i}\left(\lambda^{(i)}F_{x_{i}}  \right)
\cup_{j}\left(\mu^{(j)}F_{y_{j}} \right).
\]
}\label{fig: drawing of Cohen-Macaulay C^* fixed subscheme}
\end{figure}
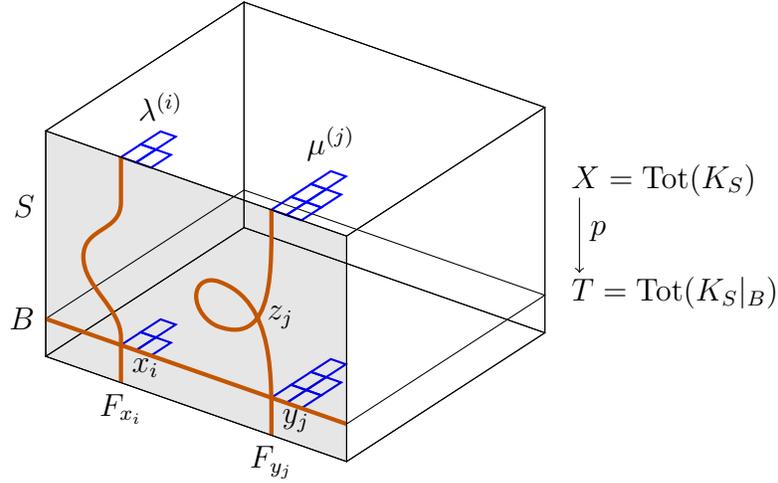

Crucially, any effective divisor on $S$ in class $[B+dF]$ is a comb
curve, i.e.~the scheme theoretic union of $B$ (our chosen section) and
some (possibly thickened) fibres with total multiplicity $d$. This is
proved in Lemma \ref{lem: Sym(B) = Hilb(S)}. The main result of this
section is the following:
\begin{theorem}\label{thm: C* invariant curves are partition thickened
comb curves with embedded points} If a subscheme $Z\subset X$ in the
class $[Z]=B+dF$ is $\CC^{*}$-invariant, then it is a partition
thickened comb curve with points. That is
\[
\Hilb^{d,n}(X)^{\CC^{*}}\subset \Hilb^{d,n}_{\PCP }(X)\subset \Hilb^{d,n}(X).
\]
Moreover, $\CC^{*}$ acts
on $\Hilb^{d,n}_{\PCP}(X)$ and there
exists a constructible morphism\footnote{A
constructible morphism is a map which is regular on each piece of a
decomposition of its domain into locally closed subsets. Because we
work with Euler characteristics and the Grothendieck group, we need
only work with constructible morphisms.}
\begin{equation}\label{eqn: rho : HilbPCP-->Sym(B)}
\rho_{d}: \Hilb^{d,\bullet}_{\PCP }(X) \to  \Sym^{d}(B)
\end{equation}
such that if $[Z]\in \Hilb_{\PCP}^{d,n}(X)$, where the maximal
Cohen-Macaulay subscheme of $Z$ is $B\cup_{i}\left(\lambda^{(i)}F_{x_{i}}
\right)$, then
\[
\rho_{d}([Z]) = \sum_{i} |\lambda^{(i)}|\,x_{i}. 
\]
\end{theorem}
\begin{proof}
We have to prove the following: Let $Z \subset X$ be a $\CC^*$-fixed 
subscheme in the class $[Z]=B+dF$, then the underlying Cohen-Macaulay 
support curve $C$ is a partition thickened comb curve. Let 
$I_C \subset \O_X$ be the ideal sheaf defining $C$. Pushing forward 
along the projection $\pi : X=K_S \rightarrow S$ and using the 
decomposition into $\CC^*$-weight spaces shows that there exist 
ideal sheaves 
\[
I_0 \subset \cdots \subset I_{l-1} \subsetneq \O_S
\]
such that
\[
\pi_* I_C = \bigoplus_{i=0}^{l-1} I_i \otimes K_{S}^{-i}.
\]
This is essentially proved in \cite[Sect.~4]{Kool-Thomas2016} (albeit
in the PT rather than the DT setting). Each $I_i$ defines a closed
subscheme $C_i \subset S$ satisfying
\begin{align*}
&S \supsetneq C_0 \supset \cdots \supset C_{l-1}, \\
&\sum_{i=0}^{l-1} [C_i] = B + dF \in H_2(S).
\end{align*}
Therefore each $C_i$ has dimension $\leq 1$. In fact each $C_i \subset S$ 
is a Cohen-Macaulay curve, or else $C$ has embedded points. Since a 
Cohen-Macaulay curve on a surface is Gorenstein, each $C_i$ is an 
effective divisor.

By the nesting condition and Lemma \ref{lem: Sym(B) = Hilb(S)}, we 
deduce
\begin{align*}
C_0 &= B + \sum_{i=1}^{n} \lambda_{1}^{(i)} F_{x_i}, \\
C_1 &=  \sum_{i=1}^{n} \lambda_{2}^{(i)} F_{x_i}, \\
&\,\, \, \,\vdots \\
C_{l-1} &= \sum_{i=1}^{n} \lambda_{l}^{(i)} F_{x_i},
\end{align*}
for some distinct points $x_1, \ldots, x_n \in B$ and 
$\lambda_{1}^{(i)} \geq \cdots \geq \lambda_{l}^{(i)}$. 
This proves that the $\CC^*$-fixed locus lies inside the PCP locus.

Since the $\CC^*$-invariant Cohen-Macaulay curves just described are
exactly the support curves of PCP curves, it follows that the PCP
locus is $\CC^*$-invariant. Finally, since the assignment $Z\mapsto
Z_{\CM}$ which takes a 1-dimensional subscheme to its maximal
Cohen-Macaulay subscheme defines a constructible morphism $\Hilb (X)
\to \Hilb (X)$, its restriction to $\Hilb^{d,\bullet }_{\PCP}(X)$ is
also constructible and thus gives the constructible morphism
$\rho_{d}$.
\end{proof}

\presectionspace
\section{Push-forward to the symmetric product} \label{sec: pushforward to sym prod}

From the $\CC^{*}$-equivariant inclusions in Theorem~\ref{thm: C*
invariant curves are partition thickened comb curves with embedded
points} and $\CC^{*}$-localization of Euler characteristic, we have
\[
\DThat (X) = \int_{\Hilb^{\bullet ,\bullet}(X)} 1 \, de  =
\int_{\Hilb^{\bullet ,\bullet}(X)^{\CC^{*}}} 1 \, de  = \int_{\Hilb^{\bullet
,\bullet}_{\PCP }(X)} 1 \, de  .
\]
We compute these Euler characteristics by pushing forward along the
map $\rho_{d}$ constructed in Theorem~\ref{thm: C* invariant curves are partition thickened
comb curves with embedded points}. That is we use

$$
\int_{\Hilb^{d,\bullet}_{\PCP }(X)} 1 \, de = \int_{\Sym^d(B)} (\rho_{d})_{*}(1) \, de,
$$
where $(\rho_{d})_{*}(1)$ is the $\ZZ (\!(p)\!)$-valued constructible
function on $\Sym^d(B)$ given by pushing forward the Euler
characteristic measure \cite{MacPherson-Annals74}. We denote
$(\rho_{d})_{*}(1)$ by $f_{d}$ so by definition, the value of $f_{d} $
at a point $\bolda \boldx =\sum_{i}a_{i}x_{i} \in \Sym^d(B)$ is
$$
f_d(\bolda \boldx ) = \int_{\rho_{d}^{-1}(\bolda \boldx )} 1 \, de.
$$

We will show that $f_{d}$ has some nice multiplicative properties. Let
$B^{\sing}\subset B$ be the points over which the fibers of $S\to B$
are singular. Note that $\# B^{\sing}=e(S)$. Let
$B^{\sm}=B-B^{\sing}$.

\begin{proposition}\label{prop: fd = F1*F2*G*H}
Let $x_{1},\dotsc ,x_{n}\in B^{\sm}$ and $y_{1},\dotsc ,y_{m}\in
B^{\sing}$ and let $a_{1},\dotsc ,a_{n},b_{1},\dotsc ,b_{m}$ be
positive integers summing to $d$. Let $\bolda \boldx$ and $\boldb
\boldy$ denote $\sum_{i}a_{i}x_{i}$ and $\sum_{j}b_{j}y_{j}$
respectively. Then there exist $F_{1} \in p^{\frac{1}{2}} \ZZ [[p]]$,
$F_{2}\in \ZZ [[p]]$, and $g,h:\NN \to \ZZ (\!(p)\!)$ such that
\[
f_{d}(\bolda \boldx +\boldb \boldy ) = F_{1}^{e(B)}\cdot
F_{2}^{e(S)}\cdot G(\bolda \boldx )\cdot H(\boldb \boldy ),
\]
where
\[
G(\bolda \boldx ) = \prod_{i=1}^{n}g(a_{i}), \quad H(\boldb \boldy ) =
\prod_{j=1}^{m}h(b_{j}).
\]
\end{proposition}
This proposition follows from Proposition~\ref{prop: formula for fd in
terms of the normalized vertex} which will be stated and proved in the
next section.

\begin{corollary}\label{cor: DThat = F1^{e(B)}F2^{e(S)}(sum g q^{a})^{e(B)-e(S)}...}
\[
\DThat (X) = F_{1}^{e(B)}\cdot F_{2}^{e(S)}\cdot
\left(\sum_{a=0}^{\infty}g(a) q^{a} \right)^{e(B)-e(S)}\cdot
\left(\sum_{b=0}^{\infty}h(b) q^{b} \right)^{e(S)},
\]
where we have set $g(0)=h(0)=1$.
\end{corollary}
\begin{proof}
We apply Proposition~\ref{prop: fd = F1*F2*G*H} to the computation of
$\DThat (X)$ as follows
\begin{align*}
\DThat (X)&= \int_{\Hilb^{\bullet ,\bullet}_{\PCP}(X)} 1\, de\\
&=\int_{\Sym^{\bullet}(B)}f_{\bullet } \, de \\
&= F_{1}^{e(B)}\cdot F_{2}^{e(S)}\cdot  \int_{\Sym^{\bullet}(B^{\sm})}
G\,de \cdot  \int_{\Sym^{\bullet}(B^{\sing})}
H\,de .
\end{align*}
Applying Lemma~\ref{lem: formula for euler char of sym products} to
this last equation yields the corollary. 
\end{proof}

To prove Proposition~\ref{prop: fd = F1*F2*G*H} and explicitly
compute $F_{1}$, $F_{2}$, $g$, and $h$, we need a good understanding
of the strata $\rho_{d}^{-1} (\bolda \boldx +\boldb \boldy )\subset
\Hilb^{d,\bullet}_{\PCP}(X)$.

For any
\begin{align*}
\boldx =(x_{1},\dotsc ,x_{n}),&\quad \boldy =(y_{1},\dotsc
,y_{m}),\\
\boldlambda =(\lambda^{(1)},\dotsc ,\lambda^{(n)}),&\quad \boldmu =(\mu^{(1)},\dotsc ,\mu^{(m)}),
\end{align*}
we define an associated Cohen-Macaulay curve
\[
C_{\boldx, \boldy, \boldlambda, \boldmu} = B \bigcup_{i=1}^{n}\left(\lambda^{(i)}F_{x_{i}} \right)
\bigcup_{j=1}^{m}\left(\mu^{(j)}F_{y_{j}} \right).
\]
From Theorem~\ref{thm: C* invariant curves are
partition thickened comb curves with embedded points} we obtain the
following decomposition of the fibers of $\rho_{d}$ in $K_0(\Var_{\CC})(\!(p)\!)$:
\begin{align} 
\begin{split} \label{eqn: components of fibers of rho}
\rho_{d}^{-1}(\bolda \boldx +\boldb \boldy ) &=  \sum_{\boldlambda
\vdash \bolda}\, \sum_{\boldmu \vdash \boldb}  p^{\chi(\O_{C_{\boldx, \boldy, \boldlambda, \boldmu}})} \Hilb^{\bullet}(X,C_{\boldx, \boldy, \boldlambda, \boldmu}) \\
&=  \sum_{\boldlambda
\vdash \bolda}\, \sum_{\boldmu \vdash \boldb}  p^{\chi(\O_{C_{\boldx, \boldy, \boldlambda, \boldmu}})} \Quot_X^{\bullet}(I_{C_{\boldx, \boldy, \boldlambda, \boldmu}}),
\end{split}
\end{align}
where the second equality follows from Lemma \ref{HilbtoQuot}. Here 
\begin{align*}
\bolda =(a_{1},\dotsc ,a_{n}),&\quad \boldb =(b_{1},\dotsc
,b_{m})
\end{align*}
and the meaning of $\boldlambda \vdash \bolda$ and $\boldmu \vdash
\boldb$ is that $\lambda^{(i)}\vdash a_{i}$ and $\mu^{(j)}\vdash
b_{j}$ for all $i$ and $j$. For later use, we state the following:

\begin{lemma}\label{lem: chi(C)=chi(B) -sum lamba1 - sum mu1}
Let 
\[
C_{\boldx, \boldy, \boldlambda, \boldmu}:=B\cup_{i}\left(\lambda^{(i)}F_{x_{i}} \right)\cup_{j}\left(\mu^{(j)}F_{y_{j}} \right),
\]
then
\[
\chi (\O_{C_{\boldx, \boldy, \boldlambda, \boldmu}}) = \chi (\O_{B}) -\sum_{i=1}^{n}\lambda^{(i)}_{1}
-\sum_{j=1}^{m}\mu^{(j)}_{1}. 
\]
\end{lemma}
\begin{proof}
Since $\lambda^{(i)}F_{x_{i}}=p^{-1}(Z_{\lambda^{(i)}})$ and $p$ is an
elliptic fibration, $\chi (\O_{F_{x_{i}}})=0$ and similarly we have
$\chi (\O_{F_{y_{j}}})=0$. Note that $B\cap \lambda^{(i)}F_{x_{i}}$
and $B\cap \mu^{(j)}F_{y_{j}}$ are zero dimensional subschemes of
length $\lambda^{(i)}_{1}$ and $\mu^{(j)}_{1} $ respectively
(c.f. equation~\eqref{eqn: monomial ideal I given by a partition
lambda}). The lemma then follows from the exact sequence
\[
0\to \O_{C} \to
\O_{B}\oplus_{i}\O_{\lambda^{(i)}F_{x_{i}}}\oplus_{j}\O_{\mu^{(j)}F_{y_{j}}}
\to \oplus_{i}\O_{B\cap \lambda^{(i)}F_{x_{i}}}\oplus_{j}\O_{B\cap \mu^{(j)}F_{y_{j}}}  \to 0.
\]
\end{proof}

In the next section, we will see that the Euler characteristic of the Quot scheme
$\Quot_X^{\bullet}(I_{C_{\boldx, \boldy, \boldlambda, \boldmu}})$ does not depend on the exact location of
the points $x_i \in B^{\sm}$ and $y_j \in B^{\sing}$, but only on
their number $n$ and $m$ and the partitions $\lambda^{(i)}$ and
$\mu^{(j)}$.

\presectionspace
\section{Stratifying according to embedded points} \label{sec: restriction
to formal nghds}

In the previous two sections we reduced our consideration to the
strata $\Quot_X^{\bullet}(I_{C_{\boldx, \boldy, \boldlambda, \boldmu}})$
of $ \Hilb^{d,\bullet}_{\PCP }(X)$ which parameterize subschemes $Z$ whose maximal
Cohen-Macaulay subscheme $Z_{\CM} \subset Z$ is the partition
thickened comb curve 
\[
C_{\boldx, \boldy, \boldlambda, \boldmu}:=B\cup_{i}(\lambda^{(i)}F_{x_{i}})\cup_{j}(\mu^{(j)}F_{y_{j}}).
\]

In this section, we introduce a further stratification of $\Quot_X^{\bullet}(I_{C_{\boldx, \boldy, \boldlambda, \boldmu}})$ 
by keeping track of the support of the quotients with respect to the geometry of the underlying reduced curve
$B\cup_{i} F_{x_{i}} \cup_{j} F_{y_{j}}$. This allows us to write $\Quot_X^{\bullet}(I_{C_{\boldx, \boldy, \boldlambda, \boldmu}})$ as a product 
of ``local'' Hilbert schemes (in $K_{0}(\Var_{\CC})(\!(p)\!)$). We then
use this product to compute its Euler characteristic. The main result of this section is
Proposition~\ref{prop: formula for fd in terms of the normalized
vertex}.

\subsection{Stratification of $X$}

Given a Cohen-Macaulay  curve $C_{\boldx, \boldy, \boldlambda, \boldmu}$, 
its reduced support is given by $B\cup_{i}F_{x_{i}}\cup_{j}F_{y_{j}}$
which is a nodal curve with nodes at $(x_{1},\dotsc ,x_{n})$,
$(y_{1},\dotsc ,y_{m})$ and $(z_{1},\dotsc ,z_{m})$ where $z_{j}$ is
the node of the nodal fiber $F_{y_{j}}$ (see Figure~\ref{fig: drawing
of Cohen-Macaulay C^* fixed subscheme}).
Consider the following associated chain of closed subsets of $X$:
\[
\cup_i \{x_i\} \cup_j \{y_j,z_j\} \subset B\cup_{i}F_{x_{i}}\cup_{j}F_{y_{j}} \subset X.
\]
This gives the following stratification of $X$ by locally closed subsets:
\begin{itemize}
\item $\cup_i \{x_i\} \cup_j \{y_j,z_j\}$,
\item $B^\circ := B \setminus \cup_i \{x_i\} \cup_j \{y_j\}$,
\item $F_{x_i}^{\circ} := F_{x_i} \setminus \{x_i\}$,
\item $F_{y_j}^{\circ} := F_{y_j} \setminus \{y_j,z_j\}$,
\item $W:=X \setminus B\cup_{i}F_{x_{i}}\cup_{j}F_{y_{j}}$.
\end{itemize}
We denote the collection of these locally closed subsets by $\Sigma_{\boldx, \boldy, \boldlambda, \boldmu}$. 

\begin{definition} \label{QuotSupp}
Let $X$ be a smooth quasi-projective variety, $\F$ a coherent sheaf on $X$, and $S \subset X$ a locally closed subset. Consider the Quot scheme $\Quot_X^n(\F)$ of quotients $\F \twoheadrightarrow Q$ on $X$, where $Q$ is zero dimensional of length $n$. We define $\Quot_X^n(\F,S)$ as the locally closed subset of quotients $\F \twoheadrightarrow Q$ for which the reduced support of $Q$ lies in $S$.
\end{definition}

We will use $\Sigma$ to provide a stratification of $\Quot_X^{\bullet}(I_{C})$ 
by locally closed subsets. For this, we need the following general result:
\begin{proposition} \label{keyprop}
Let $X$ be a smooth quasi-projective variety, $S \subset X$ a locally closed subset, 
$Z \subset X$ a closed subset, and $\F$ a coherent sheaf on $X$. Suppose 
$Z \subset S$. For any $n$, there exists a geometric constructible 
morphism\footnote{I.e.~constructible morphism which is bijective 
on $\CC$-valued points.}
\[
\Quot_X^n(\F,S) \longrightarrow \bigsqcup_{n_1+n_2=n} 
\Quot^{n_1}_X(\F,S \setminus Z) \times \Quot^{n_2}_X(\F,Z).
\]
\end{proposition}
\begin{proof}
Denote by $X_Z^{(N)}$ the $(N-1)$th order neighbourhood of $Z \subset X$. 
I.e.~let $I_Z \subset \O_X$ be the ideal defining $Z \subset X$, then 
$X_Z^{(N)} \subset X$ is the closed subscheme defined by 
$I_Z^N \subset \O_X$. Denote the inclusion $X_Z^{(N)} \subset X$ 
by $\iota$ and let $U:=X \setminus Z$.

Fix $n$. We choose $N \gg 0$ with the following property. For any 
quotient $\F \twoheadrightarrow Q$, where $Q$ is zero dimensional 
of length $\leq n$ and supported on $Z$, the \emph{scheme theoretic} 
support of $Q$ is contained in $X_Z^{(N)}$.

We now describe the map of the proposition. Given a quotient 
$\F \twoheadrightarrow Q$ in $\Quot_X^n(\F,S)$, we obtain
\begin{align*}
\F|_U &\twoheadrightarrow Q|_U \\
\iota^* \F &\twoheadrightarrow \iota^* Q.
\end{align*}
The trivial quotient $\F|_{X \setminus (\Supp Q \cap U)} \twoheadrightarrow 0$ 
and $\F|_U \twoheadrightarrow Q|_U$ glue on the overlap, so we 
obtain an element of $\Quot^{n_1}_X(\F,S \setminus Z)$ for some 
$n_1 \leq n$. Moreover, push-forward along a closed embedding 
is exact, so we obtain an element of $\Quot^{n-n_1}_X(\F,Z)$ 
as follows
\[
\F \twoheadrightarrow \iota_* \iota^* \F \twoheadrightarrow \iota_* \iota^* Q.
\]

We have to show that this constructible morphism is a bijection on 
$\CC$-valued points. We start with injectivity. Suppose 
$\F \twoheadrightarrow Q_i$, for $i=1,2$, map to the same element. 
Then they agree on $U = X \setminus Z$. So it suffices to show that 
they agree on $X \setminus (\Supp Q_1 \cap U) = X \setminus (\Supp 
Q_2 \cap U)$. By hypothesis, we know that there exists an isomorphism
\[
\begin{diagram}[h=0.9cm]
\iota^* \F       &\rTo     & \iota^* Q_1 \\
\dTo_{=}    &         &\dTo_{\cong }\\
\iota^* \F   &\rTo     & \iota^* Q_2 \\
\end{diagram}
\]
It suffices to show that pushing forward $\iota^* \F \twoheadrightarrow 
\iota^* Q_i$ to $X \setminus (\Supp Q_1 \cap U)$ and composing with 
$\F|_{X \setminus (\Supp Q_1 \cap U)} \twoheadrightarrow 
\iota_* \iota^* \F|_{X \setminus (\Supp Q_1 \cap U)}$ gives back 
$\F|_{X \setminus (\Supp Q_1 \cap U)} \twoheadrightarrow 
Q_i|_{X \setminus (\Supp Q_1 \cap U)}$. This can be checked on an 
open affine cover.

Suppose $R$ is a commutative ring (corresponding to an open affine 
subset of $X \setminus (\Supp Q_1 \cap U)$) and $I \subset R$ an ideal 
(corresponding to $Z$). Let $M \twoheadrightarrow Q$ be a quotient of 
finitely generated $R$-modules with $Q$ zero dimensional (corresponding 
to either of $\F|_{X \setminus (\Supp Q_1 \cap U)} \twoheadrightarrow 
Q_i|_{X \setminus (\Supp Q_1 \cap U)}$). By our choice of $N$, we have
$$
I^N \subset \Ann(Q) \subset I,
$$
where $\Ann(Q) \subset R$ denotes the annihilator ideal of $Q$. 
Consider the composition
\begin{equation} \label{compos}
M \twoheadrightarrow M \otimes_R R / I^N \twoheadrightarrow Q 
\otimes_R R / I^N,
\end{equation}
viewed as a morphism of $R$-modules, where the middle and third 
modules are $R$-modules via the map $R \rightarrow R / I^N$. Note 
that 
\[
Q \otimes_R R / I^N \cong Q / I^N Q \cong Q
\] 
since $I^N \subset \Ann(Q)$. Composing \eqref{compos} with this 
isomorphism gives back the original quotient $M \twoheadrightarrow Q$.

For surjectivity, take two quotients $\F \twoheadrightarrow 
Q_i$, with $Q_1$ of length $n_1$ supported on $S \setminus Z$ and $Q_2$ 
of length $n-n_1$ supported on $Z$. Then $\F|_U \twoheadrightarrow 
Q_1|_U$ and $\F|_{X \setminus \Supp Q_1} \twoheadrightarrow 
Q_2|_{X \setminus \Supp Q_1}$ agree on the overlap. 
They glue to the required quotient $\F \twoheadrightarrow Q$, 
with $Q$ of length $n$ supported on $S$.
\end{proof}

\begin{lemma}\label{lem: Sigma = product of local Hilbert schemes}
The following equation holds in $K_{0}(\Var_{\CC})(\!(p)\!)$
\begin{align*}
\Quot_X^{\bullet}(I_{C_{\boldx, \boldy, \boldlambda, \boldmu}}) = &\prod_{S \in \Sigma_{\boldx, \boldy, \boldlambda, \boldmu}} \Quot_X^{\bullet}(I_{C_{\boldx, \boldy, \boldlambda, \boldmu}},S) \\
=&\Quot_X^{\bullet}(I_{C_{\boldx, \boldy, \boldlambda, \boldmu}},W) \cdot
\Quot_X^{\bullet}(I_{C_{\boldx, \boldy, \boldlambda, \boldmu}},B^\circ) \\
&\cdot \prod_{i=1}^{n}
\Quot_X^{\bullet}(I_{C_{\boldx, \boldy, \boldlambda, \boldmu}},\{x_i\}) \cdot
\Quot_X^{\bullet}(I_{C_{\boldx, \boldy, \boldlambda, \boldmu}},F_{x_i}^\circ) \\
&\cdot \prod_{j=1}^{m}
\Quot_X^{\bullet}(I_{C_{\boldx, \boldy, \boldlambda, \boldmu}},\{y_j\}) \cdot
\Quot_X^{\bullet}(I_{C_{\boldx, \boldy, \boldlambda, \boldmu}},\{z_j\}) \\
&\qquad \cdot \Quot_X^{\bullet}(I_{C_{\boldx, \boldy, \boldlambda, \boldmu}},F^\circ_{y_j}).
\end{align*}
\end{lemma}

\begin{proof}
Let $C:= C_{\boldx, \boldy, \boldlambda, \boldmu}$. First apply the previous proposition to $\Quot_X^{\bullet}(I_C)$ 
with $S = X$ and $Z = C_{\red}$ in order to obtain
\begin{align*}
\Quot_X^{\bullet}(I_C)  =&\Quot_X^{\bullet}(I_C,W) \cdot 
\Quot_X^{\bullet}(I_C,C_{\red}).
\end{align*}
Next, apply the previous proposition to $\Quot_X^{\bullet}(I_C,C_{\red})$ 
with $S = C_{\red}$ and $Z = F_{x_1}$ in order to obtain
\begin{align*}
\Quot_X^{\bullet}(I_C)  =&\Quot_X^{\bullet}(I_C,W) \cdot 
\Quot_X^{\bullet}(I_C,C_{\red} \setminus F_{x_1}) \cdot \Quot_X^{\bullet}(I_C,F_{x_1}).
\end{align*}
Repeating this procedure, taking $S = C_{\red} \setminus F_{x_1}$ and 
$Z = F_{x_2}$ et cetera, we obtain
\begin{align*}
\Quot_X^{\bullet}(I_C) =&\Quot_X^{\bullet}(I_C,W) \cdot
\Quot_X^{\bullet}(I_C,B^\circ) \\
&\cdot \prod_{i=1}^{n}
\Quot_X^{\bullet}(I_C,F_{x_i}) \\
&\cdot \prod_{j=1}^{m}
\Quot_X^{\bullet}(I_C,F_{y_j}).
\end{align*}
Next, apply the previous proposition to $\Quot_X^{\bullet}(I_C,F_{x_1})$, 
$S = F_{x_1}$, and $Z = \{x_1\}$. Then
\begin{align*}
\Quot_X^{\bullet}(I_C) =&\Quot_X^{\bullet}(I_C,W) \cdot
\Quot_X^{\bullet}(I_C,B^\circ) \\
&\cdot \Quot_X^\bullet(I_C,F_{x_1}^\circ) \cdot \Quot_X^\bullet(I_C,\{x_1\}) \\
&\cdot \prod_{i=2}^{n}
\Quot_X^{\bullet}(I_C,F_{x_i}) \\
&\cdot \prod_{j=1}^{m}
\Quot_X^{\bullet}(I_C,F_{y_j}).
\end{align*}
Repeating for all points $x_2, \ldots, x_n, y_1, \ldots, y_m, z_1, \ldots, z_m$, 
the required identity follows.
\end{proof}

\subsection{Reduction to Quot schemes of $\CC^3$}\label{subsec: Formal coordinates and
reduction to Hilbert scheme on C3}\SubSecSpace

Let $\lambda ,\mu ,\nu$ be integer partitions which we also regard as
subsets in $(\ZZ_{\geq 0})^{2}$ by their diagram as in
\cite{Bryan-Kool-Young}. Consider the subscheme
\[
C_{\lambda \mu \nu}=C_{\lambda \emptyset \emptyset}\cup
C_{\emptyset \mu \emptyset}\cup C_{\emptyset \emptyset \nu } \subset
\CC^{3}=\Spec \CC [r,s,t]
\]
defined by the monomial ideal
\[
I_{\lambda \mu \nu} = I_{\lambda \emptyset \emptyset}\cap I_{\emptyset \mu \emptyset}\cap I_{\emptyset \emptyset \nu},
\]
where 
\begin{align*}
r^{\rho}s^{\sigma}t^{\tau}\in I_{\lambda \emptyset \emptyset} & \iff
(\sigma ,\tau )\notin \lambda ,\\
r^{\rho}s^{\sigma}t^{\tau}\in I_{\emptyset \mu  \emptyset} & \iff
(\tau,\rho  )\notin \mu ,\\
r^{\rho}s^{\sigma}t^{\tau}\in I_{\emptyset \emptyset \nu } & \iff
(\rho ,\sigma  )\notin \nu . 
\end{align*}

Consider the Quot scheme $\Quot^{n}_{\CC^3}(I_{C_{\lambda \mu \nu}})$. Inside, we have the closed subset of quotients supported set theoretically at the origin (Definition \ref{QuotSupp}), which we denote by 
\begin{equation} \label{defHilblambdamunu}
\Quot^{n}(\lambda ,\mu ,\nu) := \Quot^{n}_{\CC^3}(I_{C_{\lambda \mu \nu}},\{0\}).
\end{equation}
The kernel of such a quotient $I_{C_{\lambda \mu \nu}} \twoheadrightarrow Q$ is the ideal sheaf of a one dimensional scheme $Z$ with underlying Cohen-Macaulay curve $C_{\lambda \mu \nu}$ and its embedded points supported set theoretically at the origin. We emphasize that $Z$ need not be monomial.

We note that the permutations
$(r,s,t)\mapsto (t,r,s)$ and $(r,s,t)\mapsto (s,r,t)$ induce the
isomorphisms
\[
\Quot^{n}(\lambda ,\mu ,\nu )\cong \Quot^{n}(\nu,\lambda,\mu   ),\quad
\Quot^{n}(\lambda ,\mu ,\nu )\cong \Quot^{n}( \mu ',\lambda ',\nu ' ),
\]
where $\lambda' = \{(i,j) : (j,i) \in \lambda\}$ denotes conjugate partition.

We define $\sfVtilde_{\lambda
\mu \nu}\in \ZZ [[p]]$ by
\[
\sfVtilde_{\lambda \mu \nu} = e\left(\Quot^{\bullet}(\lambda ,\mu ,\nu ) \right)
\]
and note the symmetries 
\[
\sfVtilde_{\lambda \mu \nu}=\sfVtilde_{ \nu \lambda \mu
}=\sfVtilde_{\mu '\lambda '\nu '}.
\]

Recall that $S\subset X$ is the elliptic surface and $T=\Tot
(K_{S}|_{B}).$ For any point $p\in B$, let $R_{p} =\Tot
(K_{S}|_{F_{p}})$. We choose local formal coordinates at $p = x_i$ or $p=y_j$ such that 
\[
R_{p} = \{r=0 \},\quad S=\{s=0 \}, \quad T=\{t=0 \}
\]
and when $p$ is $z_{j}$ 
\[
R_{p}=\{rt=0 \}, \quad S=\{s=0 \}.
\]

Note that at $x_{i}$ or $y_{j}$, the curve $B$ is given by $\{s=t=0
\}$ and the fiber $F_{x_{i}}$ or $F_{y_{j}}$ is given by $\{s=r=0
\}$. At the point $z_{j}$, the fiber is a nodal curve and is given by
$\{s=rt=0 \}$.

Let $C_{\boldx, \boldy, \boldlambda, \boldmu}$ be a Cohen-Macaulay curve as in the beginning of this section and consider the Quot schemes of Lemma \ref{lem: Sigma = product of local Hilbert schemes}. 

\begin{lemma} \label{Quotboxes}
We have the following equalities in $K_0(\Var_{\CC})$
\begin{align*}
\Quot^{n}_X\left( I_{C_{\boldx, \boldy, \boldlambda, \boldmu}}, \{x_i\} \right)&=
\Quot^{n}\left(\bx ,\emptyset ,\lambda^{(i)} \right),\\
\Quot^{n}_X\left(I_{C_{\boldx, \boldy, \boldlambda, \boldmu}}, \{y_{j}\} \right)&=
\Quot^{n}\left(\bx ,\emptyset ,\mu^{(j)} \right),\\
\Quot^{n}_X\left(I_{C_{\boldx, \boldy, \boldlambda, \boldmu}}, \{z_{j}\} \right)&=
\Quot^{n}\left(\mujprime ,\emptyset ,\mu^{(j)} \right).
\end{align*}
\end{lemma}
Here $\bx$ is the unique partition of size 1 (whose diagram is a
single box) and $\emptyset$ is the empty partition. Recall that the
kernels of the quotients on the RHS do not need to be monomial ideals: although their underlying maximal Cohen-Macaulay
subscheme is monomial, they may have arbitrary embedded points at the
origin.
\begin{proof}
We prove the first equality; the others follow similarly. Let $C:=C_{\boldx, \boldy, \boldlambda, \boldmu}$. Take $N \gg 0$ such that any for any quotient $I_C \twoheadrightarrow Q$, where $Q$ is set theoretically supported at $\{x_i\}$ and of length $n$, the scheme theoretic support of $Q$ lies inside the $(N-1)$th order infinitesimal neighbourhood $X_{\{x_i\}}^{(N)}$ of $\{x_i\}$. (Recall the definition of the closed subscheme $X_{Z}^{(N)}$ from the proof of Proposition \ref{keyprop}.) Let $\iota : X_{\{x_i\}}^{(N)} \hookrightarrow X$ be the closed embedding. We have already seen in the proof of Proposition \ref{keyprop} that first restricting a quotient $I_C \twoheadrightarrow Q$ to $\iota^* I_C \twoheadrightarrow \iota^*Q$ and then pushing forward to $X$ looses no information. More precisely:
$$
I_C \twoheadrightarrow \iota_* \iota^* I_C \twoheadrightarrow \iota_* \iota^*Q
$$
is isomorphic to the original quotient $I_C \twoheadrightarrow Q$. Therefore, we have a geometric constructible morphism
\begin{equation} \label{geombij1}
\Quot^{n}_X\left( I_{C}, \{x_i\} \right) \rightarrow \Quot^{n}_{X_{\{x_i\}}^{(N)}}\left( \iota^* I_{C}, \{x_i\} \right). 
\end{equation}
Similarly, we have a geometric constructible morphism
\begin{equation} \label{geombij2}
\Quot^{n}_{\CC^3}\left( I_{C_{\bx \emptyset \lambda^{(i)}}},  \{0\} \right) \rightarrow \Quot^{n}_{(\CC^3)_{\{0\}}^{(M)}}\left( I_{C_{\bx \emptyset \lambda^{(i)}}},  \{0\} \right),
\end{equation}
where $(\CC^3)_{\{0\}}^{(M)}$ is a large enough infinitesimal neighbourhood of $0 \in \CC^3$ and we may take $M=N \gg 0$. Since $X$ is smooth and three dimensional, the formal completion of its local ring at $x_i$ is given by
$$
\widehat{\O}_{X,x_i} \cong \CC[[r,s,t]].
$$ 
Let $\mathfrak{m}_{x_i}$ be the maximal ideal corresponding to $x_i$. Then \cite[Cor.~10.4]{Atiyah-Macdonald}
$$
\widehat{\O}_{X,x_i} / \widehat{\mathfrak{m}}_{x_i}^N \cong \CC[[r,s,t]] / (r,s,t)^N \cong \CC[r,s,t] / (r,s,t)^N.
$$
Therefore $X_{\{x_i\}}^{(N)} \cong (\CC^3)_{\{0\}}^{(N)}$ and, by our choice of coordinates $r,s,t$, we have $\iota^* I_C \cong I_{C_{\bx \emptyset \lambda^{(i)}}}$. Hence
$$
\Quot^{n}_{X_{\{x_i\}}^{(N)}}\left( \iota^* I_{C}, \{x_i\} \right) \cong \Quot^{n}_{(\CC^3)_{\{0\}}^{(N)}}\left( I_{C_{\bx \emptyset \lambda^{(i)}}},  \{0\} \right)
$$
and the required equality follows from the geometric bijections \eqref{geombij1}, \eqref{geombij2}.
\end{proof}

A direct consequence of the above lemma and the symmetries of
$\sfVtilde$ is 
\begin{align}\label{eqn: e(Hilb(Xhat,Chat)=vertex}
e\left(\Quot^{\bullet}\left( I_{C_{\boldx, \boldy, \boldlambda, \boldmu}}, \{x_i\} \right) \right)&=
\sfVtilde_{\lambda^{(i)} \bx  \emptyset }, \\
e\left(\Quot^{\bullet}\left( I_{C_{\boldx, \boldy, \boldlambda, \boldmu}}, \{y_j\} \right) \right)&=
\sfVtilde_{ \mu^{(j)} \bx \emptyset }, \nonumber \\
e\left(\Quot^{\bullet}\left( I_{C_{\boldx, \boldy, \boldlambda, \boldmu}}, \{z_j\} \right) \right)&=
\sfVtilde_{\mu^{(j)} \mujprime \emptyset  }. \nonumber
\end{align}

We also choose formal local coordinates at all other points. For each
point in $B^{\circ}$, choose local coordinates $(r,s,t)$ such that
$T=\{t=0 \}$ and $S=\{s=0 \}$. For each point in $F^{\circ}_{x_{i}}$
or $F^{\circ}_{y_{j}}$, choose local coordinates $(r,s,t)$ such that
$S=\{s=0 \}$ and $R_{x_{i}}=\{r=0 \}$ or $R_{y_{j}}=\{r=0 \}$
respectively. For each point in $W$, choose any formal local
coordinates $(r,s,t)$.

Consider the following constructible support morphisms
\begin{align*}
\sigma_{W}:& \Quot^{\bullet}_X(I_{C_{\boldx, \boldy, \boldlambda, \boldmu}}, W) \to \Sym^{\bullet}(W),\\
\sigma_{B^{\circ }}:& \Quot^{\bullet}_X(I_{C_{\boldx, \boldy, \boldlambda, \boldmu}},B^{\circ}) \to \Sym^{\bullet}(B^{\circ })\\
\sigma_{F_{x_{i}}^{\circ }}:& \Quot^{\bullet}_X(I_{C_{\boldx, \boldy, \boldlambda, \boldmu}},F_{x_{i}}^{\circ}) \to \Sym^{\bullet}(F_{x_{i}}^{\circ }),\\
\sigma_{F_{y_{j}}^{\circ }}:& \Quot^{\bullet}_X(I_{C_{\boldx, \boldy, \boldlambda, \boldmu}},F_{y_{j}}^{\circ}) \to \Sym^{\bullet}(F_{y_{j}}^{\circ }).
\end{align*}
To each quotient $I_{C_{\boldx, \boldy, \boldlambda, \boldmu}} \twoheadrightarrow Q$ with kernel $I_Z \subset I_{C_{\boldx, \boldy, \boldlambda, \boldmu}} $, these maps assign the set theoretic support ---weighted by length--- of $I_{C_{\boldx, \boldy, \boldlambda, \boldmu}} / I_Z$, i.e.~the location and lengths of the embedded points or zero dimensional components of $Z$. 

Consider a point $p$ in $W$, $B^{\circ}$, $F_{x_{i}}^{\circ }$, or
$F_{y_{j}}^{\circ}$. Then using the formal local coordinates chosen
above, the same type of argument as in Lemma \ref{Quotboxes} gives the following equalities in $K_0(\Var_{\CC})$ 
\begin{align*}
\sigma_{W}^{-1}(np) &= \Quot^{n}\left(\emptyset ,\emptyset,\emptyset \right), \\
\sigma_{B^{\circ }}^{-1}(np) &= \Quot^{n}\left(\bx ,\emptyset ,\emptyset \right), \\
\sigma_{F_{x_{i}}^{\circ }}^{-1}(np) &= \Quot^{n}\left(\emptyset  ,\emptyset,\lambda^{(i)}  \right), \\
\sigma_{F_{y_{j}}^{\circ }}^{-1}(np) &= \Quot^{n}\left(\emptyset
,\emptyset,\mu^{(j)} \right).
\end{align*}
Using Proposition \ref{keyprop}, we see that the pre-images of the support morphisms satisfy the
following multiplicative property in $K_0(\Var_{\CC})(\!(p)\!)$ 
\[
\sigma^{-1}_{U}\left(\sum_{i}n_{i}p_{i} \right) = \prod_{i}
\sigma^{-1}_{U}\left(n_{i}p_{i} \right).
\]


Pushing forward the Euler characteristic measure along the support
maps, applying Lemma~\ref{lem: formula for euler char of sym
products}, and using the symmetries of $\sfVtilde$ we find the
following formulas

\begin{align}\label{eqn: e(Hilb(Xhat,Bcirc))=V^{e(Bcirc)}}
e\left(\Quot_X^{\bullet}(I_{C_{\boldx, \boldy, \boldlambda, \boldmu}},W) \right) & = \int_{\Sym^{\bullet}(W)}
\left(\sigma_{W} \right)_{*}(1)\, de = \left(\sfVtilde_{\emptyset \emptyset \emptyset }\right)^{e(W)},\\
e\left(\Quot_X^{\bullet}\left(I_{C_{\boldx, \boldy, \boldlambda, \boldmu}},B^{\circ} \right) \right) & = \int_{\Sym^{\bullet}(B^{\circ })}
\left(\sigma_{B^{\circ }} \right)_{*}(1)\, de = \left(\sfVtilde_{\bx \emptyset \emptyset }\right)^{e(B^{\circ })},\nonumber \\
e\left(\Quot_X^{\bullet}\left(I_{C_{\boldx, \boldy, \boldlambda, \boldmu}},F_{x_{i}}^{\circ} \right) \right) & = \int_{\Sym^{\bullet}(F_{x_{i}}^{\circ })}
\left(\sigma_{F_{x_{i}}^{\circ }} \right)_{*}(1)\, de = \left(\sfVtilde_{\lambda^{(i)} \emptyset  \emptyset }\right)^{e(F_{_{x_{i}}}^{\circ })},\nonumber \\
e\left(\Quot_X^{\bullet}\left(I_{C_{\boldx, \boldy, \boldlambda, \boldmu}},F_{y_{j}}^{\circ} \right) \right) & = \int_{\Sym^{\bullet}(F_{y_{j}}^{\circ })}
\left(\sigma_{F_{y_{j}}^{\circ }} \right)_{*}(1)\, de =
\left(\sfVtilde_{\mu^{(j)} \emptyset  \emptyset
}\right)^{e(F_{_{y_{j}}}^{\circ })}. \nonumber
\end{align}

We are now ready to prove the main result of this section:

\begin{proposition}\label{prop: formula for fd in terms of the
normalized vertex} Recall that $f_{d}=\left(\rho_{d} \right)_{*}(1)\in
\ZZ (\!(p)\!)$ is the push forward of the Euler characteristic measure by
the map $\rho_{d}$. As before, let $x_{1},\dotsc ,x_{n}\in
B^{\sm},y_{1},\dotsc ,y_{m}\in B^{\sing }$ and let $a_{1},\dotsc
,a_{n},b_{1},\dotsc ,b_{m}$ be positive integers summing to $d$. Then
\[
f_{d}\left(\bolda \boldx +\boldb \boldy  \right)  =
\left(p^{\half}\frac{\sfVtilde_{\bx \emptyset
\emptyset}}{\sfVtilde_{\emptyset \emptyset \emptyset}} \right)^{e(B)}
\cdot \sfVtilde^{e(S)}_{\emptyset \emptyset \emptyset} \cdot
\prod_{i=1}^{n}g(a_{i}) \prod_{j=1}^{m} h(b_{j}),
\]
where
\begin{align*}
g(a)& = \sum_{\lambda \vdash a} p^{-\lambda_{1}}
\frac{\sfVtilde_{\emptyset \emptyset \emptyset} \sfVtilde_{\lambda \bx 
\emptyset}}{\sfVtilde_{\bx \emptyset \emptyset} \sfVtilde_{\lambda
\emptyset \emptyset}},\\
h(b) &= \sum_{\mu \vdash b} p^{-\mu_{1}} \frac{\sfVtilde_{\mu \mu'
\emptyset} \sfVtilde_{ \mu \bx \emptyset}}{\sfVtilde_{\bx \emptyset
\emptyset} \sfVtilde_{\mu \emptyset \emptyset}}.
\end{align*}
Note that this proves Proposition~\ref{prop: fd = F1*F2*G*H} and
provides the values of the unknowns $g,h$ (as above) and
$F_{1},F_{2}$
\[
F_{1} = p^{\half}\,\frac{\sfVtilde_{\bx \emptyset
\emptyset}}{\sfVtilde_{\emptyset \emptyset \emptyset}} , \quad
F_{2} =\sfVtilde_{\emptyset \emptyset \emptyset}. 
\]
\end{proposition}

\begin{proof}
We apply, in order, equation~\eqref{eqn: components of fibers of rho},
Lemma~\ref{lem: Sigma = product of local Hilbert schemes},
equations~\eqref{eqn: e(Hilb(Xhat,Chat)=vertex} and \eqref{eqn:
e(Hilb(Xhat,Bcirc))=V^{e(Bcirc)}}, and Lemma~\ref{lem: chi(C)=chi(B)
-sum lamba1 - sum mu1} to compute
\begin{align*}
f_{d}(\bolda \boldx +\boldb \boldy )&= e\left(\rho_{d}^{-1}(\bolda
\boldx +\boldb \boldy) \right) \\
&= \sum_{\boldlambda \vdash \bolda}\sum_{\boldmu \vdash \boldb} p^{\chi (\O_{C_{\boldx, \boldy, \boldlambda, \boldmu}})} e\left(\Quot_X^{\bullet}(I_{C_{\boldx, \boldy, \boldlambda, \boldmu}}) \right)\\
&=\sum_{\boldlambda \vdash \bolda}\sum_{\boldmu
\vdash \boldb} p^{\chi (\O_{C_{\boldx, \boldy, \boldlambda, \boldmu}})} \cdot e(\Quot_X^{\bullet}(I_{C_{\boldx, \boldy, \boldlambda, \boldmu}}, W)) \cdot
e\left(\Quot_X^{\bullet}(I_{C_{\boldx, \boldy, \boldlambda, \boldmu}},B^{\circ} \right) \\
& \quad  \cdot \prod_{i=1}^{n}
e\left(\Quot_X^{\bullet}(I_{C_{\boldx, \boldy, \boldlambda, \boldmu}},\{x_{i}\}) \right) \cdot
e\left(\Quot_X^{\bullet}(I_{C_{\boldx, \boldy, \boldlambda, \boldmu}},F^{\circ}_{x_{i}}\right) \\
& \quad  \cdot   \prod_{j=1}^{m}
e\left(\Quot_X^{\bullet}(I_{C_{\boldx, \boldy, \boldlambda, \boldmu}},\{y_{j}\}) \right) \cdot
e\left(\Quot_X^{\bullet}(I_{C_{\boldx, \boldy, \boldlambda, \boldmu}},\{z_{j}\}) \right) \\
&\qquad \quad \cdot e\left(\Quot_X^{\bullet}(I_{C_{\boldx, \boldy, \boldlambda, \boldmu}},F^{\circ}_{y_{j}}
\right)\\
&=\sfVtilde_{\emptyset \emptyset
\emptyset}^{e(W)} \cdot \sfVtilde^{e(B^{\circ})}_{\bx \emptyset
\emptyset} \cdot  \sum_{\boldlambda \vdash \bolda}\sum_{\boldmu
\vdash \boldb} p^{\chi (\O_{C_{\boldx, \boldy, \boldlambda, \boldmu}})} \\
& \quad \cdot \prod_{i=1}^{n} \sfVtilde_{\lambda^{(i)}\bx \emptyset}
\cdot \sfVtilde^{e(F^{\circ}_{x_{i}})}_{\lambda^{(i)}\emptyset
\emptyset} \cdot \prod_{j=1}^{m} \sfVtilde_{\mu^{(j)}\bx \emptyset}
\cdot \sfVtilde_{\mu^{(j)} \mujprime \emptyset }\cdot
\sfVtilde^{e(F^{\circ}_{y_{j}})}_{\mu^{(j)}\emptyset \emptyset} \\
&=p^{\chi (\O_{B})}\cdot \sfVtilde_{\emptyset \emptyset
\emptyset}^{e(W)} \cdot \sfVtilde^{e(B^{\circ})}_{\bx \emptyset
\emptyset}\\
&\quad \cdot \prod_{i=1}^{n} \sum_{\lambda^{(i)}\vdash a_{i}}
p^{-\lambda_{1}^{(i)}} \sfVtilde_{\lambda^{(i)}\bx \emptyset} \cdot
\sfVtilde^{e(F^{\circ}_{x_{i}})}_{\lambda^{(i)}\emptyset \emptyset}\\
&\quad \cdot \prod_{j=1}^{m} \sum_{\mu^{(j)}\vdash b_{j}}
p^{-\mu_{1}^{(j)}} \sfVtilde_{\mu^{(j)} \bx \emptyset} \cdot
\sfVtilde_{\mu^{(j)} \mujprime \emptyset }\cdot
\sfVtilde^{e(F^{\circ}_{y_{j}})}_{\mu^{(j)}\emptyset \emptyset}.
\end{align*}
We note that $e(F_{x_{i}})=0$ and $e(F_{y_{j}})=1$ so that
$e(F^{\circ }_{x_{i}})=-1$ and $e(F^{\circ }_{y_{j}})=-1$. Also,
since $e(B^{\circ}) = e(B)-n-m$, we have
\begin{align*}
e(W)&= e(S) - e(B^{\circ}) - \sum_{i} e(F_{x_{i}}) - \sum_{j} e(F_{y_{j}})\\
&= e(S) - e(B) +n.
\end{align*}
The above equations allow us to redistribute the terms of
$f_{d}(\bolda \boldx +\boldb \boldy )$ as follows
\begin{align*}
f_{d}(\bolda \boldx +\boldb \boldy ) &= p^{\chi (\O_{B})} \cdot
\sfVtilde^{e(S)}_{\emptyset \emptyset \emptyset} \cdot
\left(\frac{\sfVtilde_{\bx \emptyset \emptyset}}{\sfVtilde_{\emptyset
\emptyset \emptyset}} \right)^{e(B)} \\
&\quad \cdot \prod_{i=1}^{n} \sum_{\lambda^{(i)} \vdash a_{i}}
p^{-\lambda_{1}^{(i)}} \cdot \frac{\sfVtilde_{\emptyset \emptyset
\emptyset} \sfVtilde_{\lambda^{(i)} \bx  \emptyset}}{\sfVtilde_{\bx
\emptyset \emptyset} \sfVtilde_{\lambda^{(i)}\emptyset \emptyset}} \\
&\quad \cdot \prod_{j=1}^{m} \sum_{\mu^{(j)} \vdash b_{j}}
p^{-\mu_{1}^{(j)}} \cdot \frac{\sfVtilde_{\mu^{(j)} \mujprime
\emptyset} \sfVtilde_{\mu^{(j)} \bx \emptyset}}{\sfVtilde_{\bx
\emptyset \emptyset} \sfVtilde_{\mu^{(j)}\emptyset \emptyset}}.
\end{align*}
Noting that $\chi (\O_{B}) = e(B)/2$, we see the above proves the proposition.
\end{proof}

\presectionspace
\section{Reduction to the topological vertex} \label{sec: reduction to the vertex}

In this section, we express $\DThat (X)$ in terms of the topological
vertex, and then use the trace formulas of \cite{Bryan-Kool-Young} to
obtain a closed formula for $\DThat (X)$.

\subsection{$\sfVtilde_{\lambda \mu \nu}$ in terms of $\sfV_{\lambda
\mu \nu}$.}\SubSecSpace 

Recall that the coefficients of the series $\sfVtilde_{\lambda \mu \nu}\in \ZZ [[p]]$ are
given by the Euler characteristics of Quot schemes on $\CC^3$
\[
\sfVtilde_{\lambda \mu \nu} =
e\left(\Quot^{\bullet}_{\CC^3}(I_{C_{\lambda \mu \nu}})\right). 
\]
We can compute the Euler characteristics using the $T=(\CC^{*})^{3}$-action on the Quot schemes induced by the $T$-action on
$\CC^{3}$. An ideal $I\subset \CC [r,s,t]$ is $T$-invariant if and
only if it is generated by monomials (such as $I_{C_{\lambda \mu \nu}}$). Moreover, there is a bijection
between monomial ideals and 3D partitions (see \S~6.3 of
\cite{Bryan-Kool-Young}) where a monomial ideal $I\subset \CC[r,s,t]$
corresponds to a 3D partition $\pi \in \left(\ZZ_{\geq 0} \right)^{3}$
by
\[
(\rho ,\sigma ,\tau )\in \pi  \quad \Longleftrightarrow\quad 
r^{\rho}s^{\sigma}t^{\tau} \notin I.
\]
The subschemes represented by points in
$\Quot^{\bullet}_{\CC^3}(I_{C_{\lambda \mu \nu}})^T$ are given by quotients $I_{C_{\lambda \mu \nu}} \twoheadrightarrow Q$ with kernel equal to a monomial ideal corresponding to a 3D-partition
asymptotic to $(\lambda \mu \nu )$, see
\cite[Defn~1]{Bryan-Kool-Young}. Consequently, 
\begin{align*}
\sfVtilde_{\lambda \mu \nu}& = e\left(\Quot^{\bullet}_{\CC^3}(I_{C_{\lambda \mu \nu}})^T\right)\\
&=\sum_{\pi} p^{n(\pi )},
\end{align*}
where the sum runs over all 3D partitions asymptotic to $(\lambda \mu
\nu )$ and $n(\pi )$ is the number of boxes in $\pi$ which are not
contained in any of the legs. We see that $\sfVtilde_{\lambda \mu
\nu}$ differs from the usual topological vertex by an overall
normalization
\[
\sfV_{\lambda \mu \nu} = p^{|\pi_{\min}|} \sfVtilde_{\lambda \mu \nu},
\]
where $\sfV_{\lambda \mu \nu }$ is the usual topological vertex
\cite[Defn~2]{Bryan-Kool-Young} and $|\pi_{\min}|$ is the renormalized
volume \cite[page~2]{Bryan-Kool-Young} of the minimal 3D partition asymptotic to $(\lambda \mu \nu )$.  

\begin{lemma}\label{lem: normalized vertex in terms of usual vertex}
The following hold 
\[
\sfV_{\lambda \emptyset \emptyset} =\sfVtilde_{\lambda \emptyset
\emptyset}, \quad \sfV_{\lambda \bx \emptyset} = p^{-\lambda_{1}}
\sfVtilde_{\lambda \bx \emptyset}, \quad \sfV_{\mu \mu' \emptyset}
=p^{-\Vert \mu \Vert^{2}} \sfVtilde_{\mu \mu' \emptyset},
\]
where $\Vert \mu \Vert^{2}:= \sum_{j=1}^{\infty} \mu_{j}^{2}$.
\end{lemma}
\begin{proof}
The renormalized volume of a 3D-partition asymptotic to $(\lambda \mu
\nu )$ is defined by
\[
|\pi | = \sum_{(\rho ,\sigma ,\tau )\in \pi} (1-\text{\# of legs of
$\pi$ containing $(\rho ,\sigma ,\tau )$}) .
\]
For $\pi_{\min}$, the minimal 3D-partition asymptotic to $(\lambda \bx
\emptyset)$, the only cubes contributing to $|\pi_{\min}|$ are those
contained in both the $\bx$-leg and the $\lambda$-leg. They intersect
exactly in the cubes corresponding to the first part of $\lambda$,
namely $\lambda_{1}$. Thus $|\pi_{\min}|=-\lambda_{1}$ in this case.

For the case of $(\lambda \emptyset \emptyset )$ every cube is in the
$\lambda$-leg and so $|\pi_{\min}|=0$. For the case of $(\mu \mu'
\emptyset )$, each cube in the intersection of the $\mu$-leg and the
$\mu '$-leg contribute $-1$ and all other cubes contribute 0. This
intersection is a stack of squares of side lengths
$\mu_{1},\mu_{2},\dotsc$ and hence
\[
|\pi_{\min}| = -\sum_{j=1}^{\infty} \mu_{j}^{2}. \qedhere
\]
\end{proof}

\subsection{Applying the trace formulas}\SubSecSpace 

Substituting the values of $F_{1}$, $F_{2}$, $g$, and $h$ from
Proposition~\ref{prop: formula for fd in terms of the normalized
vertex} into Corollary~\ref{cor: DThat = F1^{e(B)}F2^{e(S)}(sum g
q^{a})^{e(B)-e(S)}...}, and then substituting in the formulas from
Lemma~\ref{lem: normalized vertex in terms of usual vertex} we obtain
the following
\begin{align*}
\DThat (X) & = \left(p^{\half}\frac{\sfV_{\bx \emptyset
\emptyset}}{\sfV_{\emptyset \emptyset \emptyset}} \right)^{e(B)}\cdot
\sfV_{\emptyset \emptyset \emptyset }^{e(S)} \cdot
\left(\sum_{\lambda} \frac{\sfV_{\emptyset \emptyset \emptyset }
\sfV_{\lambda \bx \emptyset}}{\sfV_{\bx \emptyset
\emptyset}\sfV_{\lambda \emptyset \emptyset}}q^{|\lambda |}
\right)^{e(B)-e(S)}\\
&\quad \cdot \left(\sum_{\mu} p^{\Vert \mu \Vert^{2}} \frac{\sfV_{\mu \mu'
\emptyset} \sfV_{\mu \bx \emptyset}}{\sfV_{\bx \emptyset
\emptyset}\sfV_{\mu \emptyset \emptyset}}q^{|\mu |} \right)^{e(S)}.
\end{align*}
From the Okounkov-Reshetikhin-Vafa formula for the topological vertex
\cite{Ok-Re-Va}, \cite[eqn~(5)]{Bryan-Kool-Young} we get
\[
\sfV_{\emptyset \emptyset \emptyset} =M(p),\quad \sfV_{\bx \emptyset
\emptyset} = \frac{M(p)}{1-p},
\]
which we substitute into the above to find
\begin{align*}
\DThat (X) &= \left(\frac{1}{p^{-\half}-p^{\half}} \right)^{e(B)}
\left(\sum_{\lambda} (1-p)\frac{\sfV_{\lambda \bx
\emptyset}}{\sfV_{\lambda \emptyset \emptyset}} q^{|\lambda
|}\right)^{e(B)-e(S)} \\
&\quad \cdot \left(\sum_{\mu} (1-p)p^{\Vert \mu \Vert^{2}} \sfV_{\mu \mu'
\emptyset}\frac{\sfV_{\mu \bx \emptyset}}{\sfV_{\mu \emptyset
\emptyset}}q^{|\mu |} \right)^{e(S)}.
\end{align*}
Applying \cite[eqns~(2)\&(4)]{Bryan-Kool-Young}, we see that 
\begin{align*}
\DThat (X) &=\left(p^{-\half}-p^{\half} \right)^{-e(B)} \cdot
\left(\prod_{d=1}^{\infty} \frac{(1-q^{d})}{(1-pq^{d})(1-p^{-1}q^{d})} \right)^{e(B)-e(S)}\\
&\quad \cdot \left(M(p) \prod_{d=1}^{\infty}
\frac{M(p,q^{d})}{(1-pq^{d})(1-p^{-1}q^{d})} \right)^{e(S)} .
\end{align*}
Noting that $e(B)$ is even, the above expression is easily seen to be
equivalent to the formula for $\DThat (X)$ in Theorem~\ref{thm: main
thm -- formulas for DT and DTfiber}. Since the formula for
$\DThat_{\fiber}(X)$ was previously proven by Toda, we may now regard
the proof of Theorem~\ref{thm: main thm -- formulas for DT and
DTfiber} complete. In the next section, we will outline the proof of
the formula for $\DThat_{\fiber}(X)$ using our methods.

\presectionspace
\section{The case of $\DThat_{\fiber}(X)$}\label{sec: DThatfiber computation}

The formula for $\DThat_{\fiber}(X)$ given in Theorem~\ref{thm: main
thm -- formulas for DT and DTfiber} follows from a wall-crossing
computation of Toda \cite[Thm~6.9]{Toda-2012-Kyoto} along with the
PT/DT correspondence \cite{Bridgeland-PTDT}. However in this section
we describe how to adapt our computation of $\DThat (X)$ to the easier
case of $\DThat_{\fiber}(X)$. Our approach yields a proof which is
independent of Toda's.

\begin{definition}\label{defn: partition thickened fiber curve}
We say that a subscheme $C\subset X$ is a \textbf{partition thickened
fiber curve} if it is of the form
\[
C= \cup_{i}\left( \lambda^{(i)}F_{x_{i}} \right),
\]
where we are using the notation of Definition~\ref{defn: partition
thickened comb curve}. We say that a subscheme $Z\subset X$ is a
\textbf{partition thickened fiber curve with points (PFP)} if the
maximal Cohen-Macaulay subscheme $Z_{\CM}\subset Z$ is a partition
thickened fiber curve. We denote by
\[
\Hilb^{dF,n}_{\PFP}(X)\subset \Hilb^{dF,n}(X)
\]
the locus in the Hilbert scheme parameterizing partition thickened
fiber curves with points. 
\end{definition}

Our proof of Theorem~\ref{thm: C* invariant curves are partition thickened
comb curves with embedded points} is easily adapted to prove the
following:
\begin{theorem}\label{thm: C* fixed fiber curves are PFP}
If a subscheme $Z\subset X$ in the class $[Z]=dF$ is
$\CC^{*}$-invariant, then it is a partition thickened fiber curve with
points. That is
\[
\Hilb^{dF,n}(X)^{\CC^{*}} \subset \Hilb^{dF,n}_{\PFP}(X) \subset \Hilb^{dF,n}(X).
\]
Moreover, $\CC^{*}$ acts on $\Hilb^{dF,n}_{\PFP}(X)$ and there exists
a constructible morphism
\[
\rho^{\fiber}_{d}:\Hilb^{dF,\bullet}_{\PFP}(X) \to \Sym^{d}(B) 
\]
such that if $[Z]\in \Hilb^{dF,n}_{\PFP}(X)$, where the maximal
Cohen-Macaulay subscheme of $Z$ is
$\cup_{i}\left(\lambda^{(i)}F_{x_{i}} \right)$, then
\[
\rho_{d}^{\fiber}([Z]) = \sum_{i} |\lambda^{(i)}|x_{i}.
\]
\end{theorem}

Define the following Cohen-Macaulay curve
\[
C_{\boldx, \boldy, \boldlambda, \boldmu}^{\fiber} := \bigcup_{i=1}^{n}\left(\lambda^{(i)}F_{x_{i}} \right)
\bigcup_{j=1}^{m}\left(\mu^{(j)}F_{y_{j}} \right) . 
\]
Similar to equation \eqref{eqn: components of fibers of rho}, the previous theorem implies that the pre-images of points under the map
$\rho^{\fiber}_{d}$ break into components
\[
\left(\rho_{d}^{\fiber} \right)^{-1}(\bolda \boldx +\boldb \boldy ) =
\sum_{\boldlambda \vdash \bolda} \sum_{\boldmu \vdash
\boldb} p^{\chi(\O_{C_{\boldx ,\boldy ,\boldlambda ,\boldmu}^{\fiber}})} \Quot_X^\bullet(I_{C_{\boldx ,\boldy ,\boldlambda ,\boldmu}^{\fiber}}),
\]
where we have adopted the same notation as in Section~\ref{sec:
pushforward to sym prod}. 
Since the only nodes of the reduced support of $C_{\boldx, \boldy, \boldlambda, \boldmu}^{\fiber}$ are $z_{1},\dotsc ,z_{m}$ (adopting
the notation of Section~\ref{sec: restriction to formal nghds}) our
stratification is simpler in this case
\[
\Sigma^{\fiber}_{\boldx, \boldy, \boldlambda, \boldmu} =\left\{ \{z_{1}\},\dotsc ,
\{z_{m}\},F^{\circ}_{y_{1}},\dotsc ,
F^{\circ}_{y_{m}}, F_{x_{1}},\dotsc ,
F_{x_{n}}, W   \right\},
\]
where this time
\[
F^{\circ}_{y_{j}} = F_{y_{j}} - \{z_{j} \}, \quad W = X - \big( \cup_i F_{x_i} \cup_j F_{y_j} \big). 
\]

With virtually the same proof, we obtain the following analog of
Lemma~\ref{lem: Sigma = product of local Hilbert schemes}:
\begin{lemma}\label{lem: Sigmafiber = product of local Hilbert schemes}
The following equation holds in $K_{0}(\Var_{\CC})(\!(p)\!)$
\begin{align*}
\Quot_X^\bullet(I_{C_{\boldx ,\boldy ,\boldlambda ,\boldmu}^{\fiber}}) &=
\Quot_X^{\bullet}(I_{C_{\boldx ,\boldy ,\boldlambda ,\boldmu}^{\fiber}}, W) \cdot \prod_{i=1}^{n}
\Quot_X^{\bullet} \left( I_{C_{\boldx ,\boldy ,\boldlambda ,\boldmu}^{\fiber}},F_{x_{i}} \right) \\
&\quad \cdot \prod_{j=1}^{m} \Quot_X^{\bullet}
\left(I_{C_{\boldx ,\boldy ,\boldlambda ,\boldmu}^{\fiber}}, F^{\circ}_{y_{j}} \right)
\cdot \Quot_X^{\bullet}\left( I_{C_{\boldx ,\boldy ,\boldlambda ,\boldmu}^{\fiber}}, \{z_{j}\} \right). 
\end{align*}
\end{lemma}

We choose the same set of formal local coordinates at each point as we
did in Section~\ref{sec: restriction to formal nghds} and by the same
reasoning as in that section, we find
\begin{align*}
e\left(\Quot_X^{\bullet}(I_{C_{\boldx ,\boldy ,\boldlambda ,\boldmu}^{\fiber}},\{z_{j}\}) \right)&= \sfVtilde_{ \mu^{(j)} \mujprime  \emptyset },\\
e\left(\Quot_X^{\bullet}(I_{C_{\boldx ,\boldy ,\boldlambda ,\boldmu}^{\fiber}},W) \right) & = \sfVtilde^{e(W)}_{\emptyset \emptyset \emptyset },\\
e\left(\Quot_X^{\bullet}(I_{C_{\boldx ,\boldy ,\boldlambda ,\boldmu}^{\fiber}},F_{x_{i}}) \right)&= \sfVtilde^{e(F_{x_{i}})}_{\lambda^{(i)} \emptyset \emptyset },\\
e\left(\Quot_X^{\bullet}(I_{C_{\boldx ,\boldy ,\boldlambda ,\boldmu}^{\fiber}}, F^{\circ
}_{y_{j}} ) \right) & = \sfVtilde^{e(F^{\circ}_{y_{j}})}_{\mu^{(j)}
\emptyset \emptyset}.
\end{align*} 
Now since $\chi (\O_{C_{\boldx ,\boldy ,\boldlambda ,\boldmu}^{\fiber}})=0$, $e(F_{x_{i}})=e(F^{\circ}_{y_{j}})=0$,
and 
\[
e(W) = e(S) - m,
\]
we have that 
\begin{align*}
f^{\fiber}_{d}(\bolda \boldx +\boldb \boldy ) & :=
e\left(\left(\rho^{\fiber}_{d} \right)^{-1} (\bolda \boldx +\boldb \boldy ) \right)\\
&= \sum_{\boldlambda \vdash \bolda} \sum_{\boldmu \vdash \boldb} p^{\chi(\O_{C_{\boldx ,\boldy ,\boldlambda ,\boldmu}^{\fiber}})}
e\left( \Quot_X^\bullet(I_{C_{\boldx ,\boldy ,\boldlambda ,\boldmu}^{\fiber}})
\right) \\
&= \sum_{\boldlambda \vdash \bolda} \sum_{\boldmu \vdash \boldb}
\sfVtilde_{\emptyset \emptyset \emptyset}^{e(S)} \prod_{j=1}^{m}
\frac{\sfVtilde_{\mu^{(j)} \mujprime \emptyset }}{\sfVtilde_{\emptyset
\emptyset \emptyset}} .
\end{align*}

We may rewrite the above as 
\[
f^{\fiber}_{d}(\bolda \boldx +\boldb \boldy ) = \sfVtilde_{\emptyset
\emptyset \emptyset}^{e(S)} \prod_{i=1}^{n} g_{\fiber}(a_{i})
\prod_{j=1}^{m} h_{\fiber}(b_{j}),
\]
where
\[
g_{\fiber}(a)  = \sum_{\lambda \vdash a} 1 \quad \quad h_{\fiber}(b)
=\sum_{\mu \vdash b} \frac{\sfVtilde_{\mu \mu'
\emptyset}}{\sfVtilde_{\emptyset \emptyset \emptyset}}. 
\]
We then get the following result, analogous to Corollary~\ref{cor:
DThat = F1^{e(B)}F2^{e(S)}(sum g q^{a})^{e(B)-e(S)}...}, with a
similar proof
\[
\DThat_{\fiber}(X) = \sfVtilde^{e(S)}_{\emptyset \emptyset \emptyset}
\cdot \left(\sum_{\lambda}q^{|\lambda |} \right)^{e(B)-e(S)} \cdot
\left(\sum_{\mu} \frac{\sfVtilde_{\mu \mu'
\emptyset}}{\sfVtilde_{\emptyset \emptyset \emptyset}} q^{|\mu |}
\right)^{e(S)}.
\]
Substituting in the equations in Lemma~\ref{lem: normalized vertex in
terms of usual vertex}, using the well known generating function for
2D-partitions, and applying \cite[eqn~(1)]{Bryan-Kool-Young}
we get
\[
\DThat_{\fiber}(X) = M(p)^{e(S)} \cdot \left(\prod_{d=1}^{\infty}
(1-q^{d})^{-1} \right)^{e(B)-e(S)} \cdot  \left(\prod_{d=1}^{\infty}
(1-q^{d})^{-1} M(p,q^{d}) \right)^{e(S)} 
\]
which is easily seen to be equivalent to the formula for
$\DThat_{\fiber}(X)$ given in Theorem~\ref{thm: main thm -- formulas
for DT and DTfiber}. \qed

\presectionspace
\section{Including the Behrend function} \label{sec: Behrend}

The aim of this section is to prove Theorem~\ref{thm: DT(X) assuming
the Behrend function conjecture}, which says that up to an overall
sign, the partition functions $\DThat (X)$ and $\DT (X)$ are equal
after the simple change of variables $y=-p$. In order to do this we
will need to assume a conjecture about the Behrend function which we
formulate below for general Calabi-Yau threefolds and may be of
independent interest.

Let $Y$ be any quasi-projective Calabi-Yau
threefold.  Let $C\subset Y$ be a (not necessarily reduced)
Cohen-Macaulay curve with proper support. Assume that the
singularities of $C_{\red}$ are locally toric\footnote{This means that
formally locally $C_{\red}$ is either smooth, nodal, or the union of
the three coordinate axes. That is at $p\in C_{\red}\subset Y$ the
ideal $\widehat{I}_{C_{\red}}\subset \widehat{\mathcal{O}}_{Y,p}$ is
given by $(x_{1},x_{2})$, $(x_{1},x_{2}x_{3})$, or
$(x_{1}x_{2},x_{2}x_{3},x_{1}x_{3})$ for some isomorphism
$\widehat{\mathcal{O}}_{Y,p}\cong \CC
[[x_{1},x_{2},x_{3}]]$. }. Recall that by Definition~\ref{defn:
Hilb(U,C)} 
\[
\Hilb^{n}(Y,C) = \{Z\subset Y \text{ such that $C\subset Z$ and
$I_{C}/I_{Z}$ has finite length $n$} \}.
\]

Note that $\Hilb^{n}(Y,C)\subset \Hilb (Y)$ and let $\nu$ denote the
Behrend function on $\Hilb (Y)$. Our conjecture is the following:

\begin{conjecture}\label{conj: Behrend fnc conj}
\[
\int_{\Hilb^{n}(Y,C)} \nu \, de = (-1)^{n} \nu ([C]) \int_{\Hilb^{n}(Y,C)} \, de,
\]
where $\nu ([C])$ is the value of the Behrend function at the point $[C]\in \Hilb (Y)$.
\end{conjecture}

\begin{remark}
Conceivably, the condition that $C_{\red}$ has locally toric
singularities could be weakened, although we do not have any evidence
for this case. Our conjecture is true for $Y$ a (globally) toric
Calabi-Yau, with torus $T$, and $C \subset Y$ any $T$-fixed
Cohen-Macaulay curve. In this case both sides of our conjecture can be
computed by restricting to the $T$-fixed points, which are
isolated. The Behrend function at such a fixed point $P$ is given by $\nu(P) =
(-1)^{\dim T_P\Hilb (Y)}$, \cite[Thm.~3.4]{Behrend-Fantechi08}. The
calculation of $(-1)^{\dim T_{P}\Hilb (Y)}$ is done in \cite{MNOP1};
specifically, the LHS of the equation in \cite[Thm~2]{MNOP1} is easily
seen to be $(-1)^{\dim T_{P}\Hilb (Y)}$ while the RHS clearly obeys the
formula in our conjecture.

One could also make the much stronger conjecture that 
\[
\nu ([Z]) = (-1)^{n} \nu ([C]),
\]
for all $[Z]\in \Hilb^{n}(Y,C)$. This would of course imply our
conjecture as stated. However, we do not know whether this stronger version
holds, even in the case where $Y$ is $\CC^{3}$ and $C$ is empty. In
this case, this stronger conjecture says that the Behrend function
on $\Hilb^{n}(\CC^{3})$ is the constant function $(-1)^{n}$.
\end{remark}

\subsection{Proof of Theorem~\ref{thm: DT(X) assuming the Behrend
function conjecture}}\SubSecSpace 

The Behrend function of any scheme is invariant under
automorphisms. In particular, it is constant on the orbits of the
$\CC^{*}$ action on $\Hilb (X)$. We thus have
\[
\DT (X) = \int_{\Hilb^{\bullet ,\bullet}(X)} \nu \, de =
\int_{\Hilb^{\bullet ,\bullet}(X)^{\CC^{*}}} \nu \, de =
\int_{\Hilb_{\PCP }^{\bullet ,\bullet}(X)} \nu \, de
\]
and so
\[
\DT (X)  =\int_{\Sym^{\bullet}(B)} (\rho_{\bullet})_{*} (\nu )\, de.
\]
Let $f^{\nu}_{d} = (\rho_{d})_{*}(\nu )$ so that in the notation of
Section~\ref{sec: pushforward to sym prod}, we have
\[
f^{\nu}_{d}(\bolda \boldx +\boldb \boldy ) =
\int_{\rho^{-1}_{d}(\bolda \boldx +\boldb \boldy )} \nu \, de.
\]
Recall that for the partition function $\DT (X)$, the variable
tracking the holomorphic Euler characteristic is $y$ rather than $p$
so $f_{d}^{\nu}(\bolda \boldx +\boldb \boldy )\in \ZZ (\!(y)\!)$.

In Section \ref{sec: pushforward to sym prod}, at the level of $\CC$-valued points, we expressed $\rho^{-1}_{d}(\bolda \boldx +\boldb \boldy )$ as a disjoint union of closed subsets $\Hilb^\bullet(X,C_{\boldx ,\boldy ,\boldlambda ,\boldmu})$. We obtain
\[
f_{d}^{\nu}(\bolda \boldx +\boldb \boldy ) = \sum_{\boldlambda \vdash
\bolda} \sum_{\boldmu \vdash \boldb} \,\, y^{\chi(\O_{C_{\boldx ,\boldy ,\boldlambda ,\boldmu}})+n}\int_{\Hilb^n(X,C_{\boldx ,\boldy ,\boldlambda ,\boldmu})} \nu \, de,
\]
where $\nu$ is the Behrend function of $\Hilb^{\bullet,\bullet}(X)$ and $C_{\boldx, \boldy, \boldlambda, \boldmu}$ denotes the following Cohen-Macaulay curve 
\[
C_{\boldx, \boldy, \boldlambda, \boldmu}:=B \cup_{i}\left(\lambda^{(i)}F_{x_{i}} \right) \cup_{j}\left(\mu^{(j)}F_{y_{j}} \right).
\]
Recall that the factor $y^{\chi(\O_{C_{\boldx ,\boldy ,\boldlambda ,\boldmu}})}$ comes from the fact that $\Hilb^{\bullet,\bullet}(X)$ is indexed by $\chi(\O_Z)$ and $\Hilb^\bullet(X,C_{\boldx ,\boldy ,\boldlambda ,\boldmu})$ by the length of $I_{C_{\boldx ,\boldy ,\boldlambda ,\boldmu}} / I_Z$ (Definition \ref{defn: Hilb(U,C)}). 
We apply Conjecture~\ref{conj: Behrend fnc conj} to the above and also
use 
\[
\nu ([C_{\boldx ,\boldy ,\boldlambda ,\boldmu}]) = (-1)^{\chi (\O_{S})-\chi (\O_{C_{\boldx ,\boldy ,\boldlambda ,\boldmu}})}, 
\]
the highly non-trivial result given in Corollary~\ref{cor: value of Behrend
function at comb curves} and proved in the next section. We find
\begin{align*}
f_{d}^{\nu}(\bolda \boldx +\boldb \boldy )& = \sum_{\boldlambda \vdash
\bolda} \sum_{\boldmu \vdash \boldb} y^{\chi
(\O_{C_{\boldx ,\boldy ,\boldlambda ,\boldmu}})}\sum_{n=0}^{\infty} y^{n}(-1)^{n}\nu ([C_{\boldx ,\boldy ,\boldlambda ,\boldmu}]) \int_{\Hilb^{n}(X,C_{\boldx ,\boldy ,\boldlambda ,\boldmu})}  \, de\\
&=(-1)^{\chi (\O_{S})} \sum_{\boldlambda \vdash \bolda} \sum_{\boldmu
\vdash \boldb} \sum_{n=0}^{\infty} (-y)^{\chi (\O_{C_{\boldx ,\boldy ,\boldlambda ,\boldmu}})+n}
\int_{\Hilb^{n}(X,C_{\boldx ,\boldy ,\boldlambda ,\boldmu})} \, de.
\end{align*}
After the substitution $-y=p$, we find that the above is equivalent to 
\[
f^{\nu}_{d}(\bolda \boldx +\boldb \boldy ) = (-1)^{\chi (\O_{S})}
f_{d}(\bolda \boldx +\boldb \boldy ).
\]
It follows that 
\[
\DT (X) = (-1)^{\chi (\O_{S})} \DThat (X)
\]
after the substitution $p=-y$ as asserted by Theorem~\ref{thm: DT(X)
assuming the Behrend function conjecture}.

The case of $\DT_{\fiber}(X)$ (previously shown by Toda) proceeds
similarly except that it does not require the difficult deformation
result of the next section. Indeed, in this case, we only need to know
that the value of the Behrend function at a partition thickened fiber
curve is 1
\[
\nu \left(\left[ \cup_{i} \left(\lambda^{(i)}F_{x_{i}} \right) \right]
\right) =1.
\]
This follows from the fact that subschemes in $X$ of the form
$p^{-1}(Z)$, where $Z$ is a zero dimensional subscheme of $T$, form a
component of $\Hilb (X)$ which is isomorphic to the Hilbert scheme of
points on $T$ and hence smooth and even dimensional. While this can be
proven directly, one can also do a similar (but easier) computation as
we do in the proof of Theorem~\ref{thm: Ext computation} in
Section~\ref{sec: smoothness and deformations}.

\presectionspace
\section{Smoothness and infinitesimal deformations}\label{sec: smoothness and deformations}

In this section we show that the locus of partition thickened comb
curves lies in the non-singular locus of $\Hilb (X)$ and we compute
the dimension of $\Hilb (X)$ at those points. As a corollary, we
determine the value of the Behrend function at the points of the
Hilbert scheme corresponding to partition thickened comb curves. This
is the key technical result required in Section~\ref{sec: Behrend} to
promote our computation of $\DThat (X)$ to a computation of the
Behrend function version $\DT (X)$.

We begin by stating the three main results of this section.

\begin{theorem}\label{thm: strata of Hilb(C^2) with fixed intersection
is smooth} Let $B\subset T$ be a smooth curve contained in a smooth
surface $T$. Define $V_{l}\subset \Hilb^{d}(T)$ to be the locus of
points parameterizing zero dimensional subschemes $Z\subset T$ of length
$d$ such that the length of the scheme theoretic intersection $Z\cap
B$ is $l$. Then $V_{l}$ is locally closed and smooth of dimension
$2d-l$.
\end{theorem}

\begin{theorem}\label{thm: Ext computation}
Let $\lambda^{(1)},\dotsc ,\lambda^{(n)}$ be partitions and let
$C=B\cup_{i}\left(\lambda^{(i)}F_{x_{i}} \right)$ be a partition
thickened comb curve. The Zariski tangent space of $\Hilb (X)$ at
the point $[C]$, which is given by $\Hom (I_{C},\O_{C})\cong
\Ext^{1}_{0}(I_{C},I_{C})$, has dimension
\[
h^{0}(N_{B/T}) + \sum_{i=1}^{n}\left(2|\lambda^{(i)}| - \lambda_{1}^{(i)} \right).
\]
\end{theorem}

\begin{theorem}\label{thm: locus of comb curves is smooth}
The locus of partition thickened comb curves is contained in the
non-singular locus of $\Hilb (X)$.
\end{theorem}

\begin{corollary}\label{cor: value of Behrend function at comb curves}
The value of the Behrend function at $[C]\in \Hilb (X)$ for a
partition thickened comb curve $C =
B\cup_{i}\left(\lambda^{(i)}F_{x_{i}} \right)$ is given by
\[
\nu ([C]) = (-1)^{\chi (\O_{S})-\chi (\O_{C})}.
\]
\end{corollary}
\begin{proof}
By \cite{Behrend-micro}, the Behrend function on a smooth scheme $V$
is $(-1)^{\dim V}$ and so by Theorems~\ref{thm: locus of comb curves
is smooth} and \ref{thm: Ext computation}
\[
\nu ([C]) = (-1)^{h^{0}(N_{B/X}) }\prod_{i=1}^{n} (-1)^{\lambda_{1}^{(i)}}. 
\]
Lemma~\ref{lem: h0 of normal bundles of B in S and T} and
Lemma~\ref{lem: chi(C)=chi(B) -sum lamba1 - sum mu1} say that
\[
h^{0}(N_{B/X}) = \chi (\O_{S}) - \chi (\O_{B}), \quad \chi (\O_{C}) =
\chi(\O_{B}) - \sum_{i} \lambda^{(i)}_{1}
\]
which, when substituted into the above, prove the corollary. 
\end{proof}
The most difficult of the above results is Theorem~\ref{thm: Ext
computation} and its proof occupies the majority of this section.

Our method for computing the dimension of deformation spaces is an
adaption of Haiman's method for computing infinitesimal deformations
of zero dimensional subschemes on a surface \cite{Haiman1998}. Indeed,
the proof of Theorem~\ref{thm: strata of Hilb(C^2) with fixed
intersection is smooth} follows directly using Haiman's argument. For
Theorem~\ref{thm: Ext computation}, we use Haiman's method to study
local deformations of $C$ in the formal neighborhoods of the points
$x_{i}$, but we use the global geometry to keep track of which local
deformations extend.

\subsection{Setup for the proof of Theorem~\ref{thm: Ext
computation}.}\SubSecSpace 

For notational simplicity we first treat the case where there is a
single partition thickened fiber $F=F_{x}$ at $x\in B$, that is
\[
C=B\cup \lambda F,
\]
where $(\lambda_{1}\geq \lambda_{2}\geq \dotsb \geq \lambda_{l})$ is
an integer partition of length $l$.

Consider the divisors
\[
S,\quad R=\Tot (K_{S}|_{F}),\quad T=\Tot (K_{S}|_{B})
\]
and let $(r,s,t)$ be formal local coordinates at $x$ such that 
\[
R=\left\{r=0 \right\},\quad S=\left\{s=0 \right\},\quad T=\left\{t=0 \right\}.
\]
The formal local ring $\widehat{\O}_{X,x}\cong \CC [[r,s,t]]$ has a basis
as a $\CC$-vector space given by monomials
$\left\{r^{\rho}s^{\sigma}t^{\tau} \right\}$ for $(\rho, \sigma,
\tau)\in \left(\ZZ_{\geq 0} \right)^{3}$. We visualize these basis
vectors as unit cubes in the positive octant of $\RR^{3}$ with the
monomial $r^{\rho}s^{\sigma}t^{\tau}$ corresponding to the cube whose
corner closest to the origin is at $(\rho ,\sigma ,\tau)$.

\subsection{Exact sequences}\SubSecSpace 

The ideal sheaf $I_{C}$ has a locally free resolution of the form
\begin{equation}\label{eqn: R-->G-->I_C exact sequence}
0\to \oplus_{\beta}R_{\beta} \to \oplus_{\alpha }G_{\alpha } \to I_{C}
\to 0,
\end{equation}
where $G_{\alpha}$ (the ``generators'') and $R_{\beta}$ (the
``relations'') are of the form 
\[
\O (-\rho R-\sigma S-\tau T).
\]

Indeed, we can explicitly take the collection of $(\rho ,\sigma
,\tau)$ for $G_{\alpha}$ to be 
\[
\left\{(\lambda_{1},0,1),(\lambda_{2},1,0),(\lambda_{3},2,0),\dotsc
,(\lambda_{l},l-1,0),(0,l,0) \right\}.
\]
Note that the $\tau$ component is 1 for the first generator, and zero
for all others.

We also have the sequence
\[
0\to \O_{C} \to \O_{B}\oplus\O_{\lambda F}\to \O_{\lambda_{1}x}\to 0,
\]
where $\lambda_{1}x=B\cap \lambda F$ is the length $\lambda_{1}$
subscheme of $B$ supported at $x$.

By standard homological algebra, we have that $\Hom (I_{C},\O_{C})$ is
$H^{0}$ of the complex
\[
\Hom \left(\left[\oplus_{\beta}R_{\beta}\to
\oplus_{\alpha}G_{\alpha} \right],[\O_{B}\oplus \O_{\lambda F}\to
\O_{\lambda_{1}x}] \right).
\]
Namely, we have that $\Hom (I_{C},\O_{C})$ is given by the kernel of
the map
\[
\Hom (\oplus_{\alpha}G_{\alpha},\O_{B}\oplus \O_{\lambda
F}) \xrightarrow{\Phi_{1}\oplus \Phi_{2}}
\Hom(\oplus_{\alpha}G_{\alpha},\O_{\lambda_{1}x})\oplus \Hom
(\oplus_{\beta}R_{\beta},\O_{B}\oplus \O_{\lambda F}).
\]
This identification of $\Hom (I_{C},\O_{C})$ has a straightforward
interpretation: a homomorphism $I_{C}\to \O_C$ is determined by what
it is on each of the generators of $I_{C}$, considered as maps to
$\O_{B}$ and to $\O_{\lambda F}$. To be in the kernel of $\Phi_{1}$
just means that these maps should agree on $B\cap \lambda F$ and to be
in the kernel of $\Phi_{2}$ means that the images must obey the module
relations. We will make this combinatorially more explicit by studying
the restriction of the homomorphisms $\oplus_{\alpha}G_{\alpha}\to
\O_{B}\oplus \O_{\lambda F}$ to $\widehat{X}_{x}\cong \Spec \CC
[[r,s,t]]$.

\subsection{Combinatorics of Haiman arrows}\label{subsec: combinatorics of Haiman arrows}
\begin{figure}
\input{3Dpartition.tex}
\caption{Monomials in the local ring
$\widehat{\O}_{X,x} \cong \CC [[r,s,t]]$ are represented by
cubes. Cubes shown are the monomials in
$\widehat{\O}_{C,x}$. The grey balls are located at monomials which
generate $\hat{I}_{C,x}$ }\label{fig: B union lambda F cubes}
\end{figure}
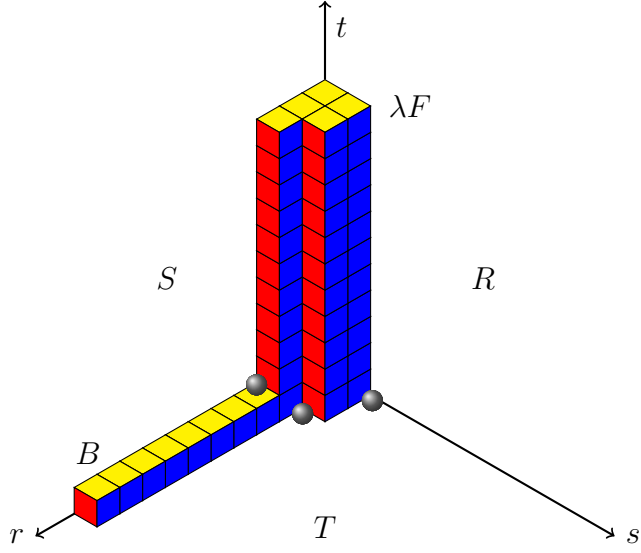

\SubSecSpace 
When restricted to the local ring $\widehat{\O}_{X,x} \cong \CC [[r,s,t]]$,
$\O_{C}$ is spanned over $\CC$ by the monomials
$r^{\rho}s^{\sigma}t^{\tau}$, where $(\rho ,\sigma ,\tau )$ are of the
form $(\rho ,0,0)$ or $(\rho ,\sigma ,\tau )_{(\rho ,\sigma )\in
\lambda}$ and $I_{C}$ is spanned by the complementary monomials. As
previously discussed, we view these monomials as cubes in the positive
octant, see Figure~\ref{fig: B union lambda F cubes}.

We call the cubes corresponding to $(\rho ,0,0)$ and $(\rho ,\sigma
,\tau )_{(\rho ,\sigma )\in \lambda}$ the $B$-cubes and $\lambda
F$-cubes respectively and the cubes in the union are called
$C$-cubes. The complement of the $C$-cubes are the $I_{C}$-cubes.

A \textbf{Haiman arrow} 
\[
(\rho ,\sigma ,\tau )\to (\rho ',\sigma ',\tau ')
\]
is a vector whose tail $(\rho ,\sigma ,\tau )$ is an $I_{C}$-cube and
whose head $(\rho ',\sigma ',\tau ')$ is a $C$-cube.

The Haiman arrows form a basis for the $\CC$-linear maps from
$\widehat{I}_{C,x}$ to $\hat{\O}_{C,x}$. We wish to determine a basis
for $\Hom (I_{C},\O_{C})$ in terms of Haiman arrows.

The generators of $I_{C}$ correspond to the cubes in the corners of
the set of $I_{C}$-cubes. They are located at $(\rho ,\sigma ,\tau )$
where $(\rho ,\sigma )$ are the corners just outside of $\lambda$ and
$\tau =0$ unless $\sigma =0$ in which case $\tau =1$ (they are
indicated by the grey balls in Figure~\ref{fig: B union lambda F
cubes}). A generator at $(\rho ,\sigma ,\tau )$ corresponds to the
image of $G_{\alpha}\to \O$ where $G_{\alpha}\cong \O (-\rho R-\sigma
S-\tau T)$. The summands of
\[
\Hom (\oplus_{\alpha}G_{\alpha},\O_{B}\oplus \O_{\lambda F}) \cong
\oplus_{\alpha} H^{0}(B,G_{\alpha}^{\vee}\otimes \O_{B})\oplus
H^{0}(F,G_{\alpha}^{\vee}\otimes \O_{\lambda F})
\]
are interpreted as follows. For $G_{\alpha}\cong \O (-\rho R-\sigma
S-\tau T)$, a homomorphism in $\Hom (G_{\alpha},\O_{B})$ or $\Hom
(G_{\alpha},\O_{\lambda F})$ is determined by a linear combination of
Haiman arrows from $(\rho ,\sigma ,\tau )$ to (respectively) some
$B$-cube or $\lambda F$-cube $(\rho ',\sigma ',\tau ').$ The location
of the head of such a Haiman arrow is determined by the order of
vanishing of the corresponding section of
$H^{0}(B,G_{\alpha}^{\vee}\otimes \O_{B})$ or
$H^{0}(F,G_{\alpha}^{\vee}\otimes \O_{\lambda F})$ --- the head will
occur at $(\rho ',\sigma ',\tau ')$ if the corresponding section is
order $r^{\rho '}s^{\sigma '}t^{\tau '}$.

We wish to determine a basis for $\Hom (I_{C},\O_{C})\cong \Ker
(\Phi_{1}\oplus \Phi_{2})$ in terms of Haiman arrows. To be in the
kernel of $\Phi_{1}$ just means that a Haiman arrow whose head is both
a $B$-cube and a $\lambda F$-cube must arise as sections of both
$H^{0}(B,G_{\alpha}^{\vee}\otimes \O_{B})$ and
$H^{0}(F,G_{\alpha}^{\vee}\otimes \O_{\lambda F})$. As for the kernel
of $\Phi_{2}$, the key observation is the following, essentially due
to Haiman \cite{Haiman1998}:
\begin{remark}\label{rem: ker of phi2 gives equations equating
translatedarrows} The equations defining the kernel of $\Phi_{2}$
equate two Haiman arrows which are obtained from one another by
translation through other Haiman arrows. Moreover, if a Haiman arrow
can be translated so that its head passes into an octant with a
negative coordinate (without its tail ever leaving the $I_{C}$-cubes)
then in must be zero.
\end{remark}

We now analyze the possible equivalence classes of Haiman arrows.

\subsection{Haiman arrows to $\lambda F$-cubes}
\SubSecSpace 
Let $G_{\alpha}\cong \O (-\rho R-\sigma S-\tau T)$ be a generating
line bundle and consider the sections
$H^{0}(F,G_{\alpha}^{\vee}\otimes \O_{\lambda F})$. A basis for this
vector space corresponds to the possible Haiman arrows $(\rho ,\sigma
,\tau )\to (\rho ',\sigma ',\tau ')$ to $\lambda F$-cubes. Since the
normal bundle of $F$ in $X$ is trivial, $\O (R)$ and $\O (S)$ are
trivial restricted to $F$. Thus
\[
G_{\alpha}^{\vee}\otimes \O_{\lambda F} \cong \O_{\lambda F}(\rho
R+\sigma S+\tau T)\cong \O_{\lambda F}(\tau x).
\]
Since $\tau$ is either 0 or 1 for the generators $G_{\alpha}$, the
Haiman arrows correspond to 
\[
H^{0}(F,\O_{\lambda F})\text{  if  }\tau =0,\quad H^{0}(F,\O_{\lambda
F}(x))\text{  if  }\tau =1.
\]
In both cases, this space has a basis of sections which in the local
coordinates are given by $\{r^{\rho '}s^{\sigma '}t^{\tau } \}_{(\rho
',\sigma ')\in \lambda}$. Note that the sections we consider above are
uniquely determined by their value on the formal neighborhood
$\hat{X}_{x}\cong \Spec \CC [[r,s,t]]$, a property which uses
crucially the fact that the genus of $F$ is 1.

We have seen that the possible Haiman arrows to $\lambda F$-cubes are
given by 
\[
(\rho ,\sigma ,\tau )\to (\rho ',\sigma ',\tau '),
\]
where $(\rho ,\sigma ,\tau )$ is a generating cube, $\tau '=\tau$ and
$(\rho ',\sigma ')\in \lambda$. In particular, the direction of the
arrows is horizontal since there is no $\tau$ component in $(\rho
-\rho ',\sigma -\sigma ',0)$.

Since all the Haiman arrows to $\lambda F$-cubes are horizontal, we
view them from above in the $(r,s)$ plane (see Figure~\ref{fig: Haiman arrows
to F-cubes}).

\begin{figure}
\begin{tikzpicture}
\draw [ultra thick,fill=yellow] (0,2)--(2,2)--(2,1)--(3,1)--(3,0)--(0,0)--(0,2);
\draw [thick,->](0,0)--(8,0);
\draw [thick,->](0,0)--(0,4);
\draw  (0,0) rectangle (1,1);
\draw  (1,0) rectangle (2,1);
\draw  (2,0) rectangle (3,1);
\draw  (0,1) rectangle (1,2);
\draw  (1,1) rectangle (2,2);
\draw [dashed] (3,1)--(8,1);
\draw [dashed] (4,1)--(4,0);
\draw [dashed] (5,1)--(5,0);
\draw [dashed] (6,1)--(6,0);
\draw [dashed] (7,1)--(7,0);
\node [left] at (0,4) {$s$};
\node [right] at (8,0) {$r$};
\draw [ultra thick,green,->] (0.55,2.5)--(0.55,1.5);
\draw [ultra thick,red,->] (0.45,2.5)--(0.45,0.5);
\draw [ultra thick,red,->] (2.5,1.5)--(1.5,0.5);
\draw [ultra thick,green,->] (3.5,0.5)--(2.5,0.5);
\draw [ultra thick,green,->] (3.5,0.5)--(1.5,1.5);
\shade [ball color=gray] (0.5,2.5) circle [radius=0.2cm];
\shade [ball color=gray] (2.5,1.5) circle [radius=0.2cm];
\shade [ball color=gray] (3.5,0.5) circle [radius=0.2cm];

\end{tikzpicture}
\caption{Examples of Haiman arrows to $\lambda F$-cubes. The green
arrows are non-zero and the red arrows are necessarily zero.}\label{fig: Haiman arrows
to F-cubes}
\end{figure}
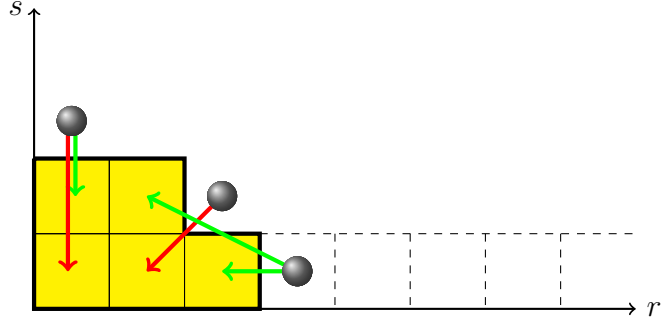

If the direction of the Haiman arrow is strictly southwest (i.e.~it
has strictly negative $\rho$ and $\sigma$ components), then by
translating (see Remark~\ref{rem: ker of phi2 gives equations equating
translatedarrows}) along the contour of $\lambda$ to the edge of the
$s$-axis, the arrow can be equated to an arrow whose head has a
negative $\rho$ component and is thus zero. There are no strictly
northeast pointing Haiman arrows, so all non-zero Haiman arrows must
be weakly northwest pointing or weakly southeast pointing. Translating
a weakly northwest pointing arrow as far to the northwest as possible,
we find that its head will either cross the $s$-axis (and hence be 0)
or it will be at the top of a column of $\lambda$ and its tail just
outside a row. Indeed, for each square in $\lambda $, there is exactly
one equivalence class of weakly northwest pointing Haiman arrows
represented by the arrow going from just outside the box's row to the
top of the box's column. Similarly, there is one equivalence class of
weakly southeast pointing Haiman arrows for each box in $\lambda$
represented by the arrow going from just outside the top of the
box's column to the end of the box's row.

The above accounts for precisely $2|\lambda |$ different equivalence
classes of Haiman arrows to $\lambda F$-cubes. However, $\lambda_{1}$
of these arrows have their head in a $B$-cube, namely the southeast
pointing arrows whose tails are just above the top of $\lambda$ and
whose head is the last square in the first row of $\lambda$. Note that
the northwest pointing Haiman arrows whose heads are in the first row
of $\lambda$ are necessarily strictly west pointing and hence
originate at the generator whose $\tau$ component is 1. Therefore the
head of these arrows also have $\tau$ component 1 and so they are not
$B$-cubes.

We thus have exactly $2|\lambda |-\lambda_{1}$ distinct equivalence
classes of Haiman arrows to $\lambda F$-cubes which are not also
arrows to $B$-cubes.

\subsection{Haiman arrows to $B$-cubes.}\SubSecSpace Any non-zero
Haiman arrow to a $B$-cube must have a tail with coordinates $(\rho
,0,1)$ or $(\rho ,1,0)$ since if not, it could be translated (see
Remark~\ref{rem: ker of phi2 gives equations equating
translatedarrows}) to an arrow whose head has negative $\tau$ or
$\sigma$ coordinates by first translating sufficiently far in the
positive $\rho$-direction and then translating the tail so that it is
just outside of the $B$-cubes. A Haiman arrow to a $B$-cube whose tail
is $(\rho ,0,1)$ or $(\rho ,1,0)$ corresponds respectively to a
section in $H^{0}(B,\O_{B}(\rho R+T))$ or $H^{0}(B,\O_{B}(\rho
R+S))$. Since
\[
\O_B (R)\cong \O_{B}(x),\quad \O_{B}(T)\cong N_{B/S},\quad \O_{B}(S)\cong N_{B/T},
\]
we see that the Haiman arrows from $(\rho ,0,1)$ or $(\rho ,1,0)$ to
$B$-cubes are given by 
\[
H^{0}(B,N_{B/S}(\rho x)) \quad \text{or}\quad
H^{0}(B,N_{B/T}(\rho x))
\]
respectively. The head of such a Haiman
arrow is $(\rho ',0,0)$ where $\rho '$ is the order of vanishing at
$x$ of the corresponding section.

\begin{lemma}\label{lem: Haiman arrows to B-cubes must vanish to the
order of the generator} Let $(\rho ,0,1)\to (\rho ',0,0)$ or $(\rho
,1,0)\to (\rho ',0,0)$ be a non-zero Haiman arrow. Then $\rho '\geq
\rho$.
\end{lemma}
\begin{proof}
Consider a Haiman arrow $(\rho ,0,1)\to (\rho ',0,0)$ with $\rho
'<\rho$. Then this arrow can be translated so that its head is a
$\lambda F$-cube, however we saw in the previous subsection that
Haiman arrows to $\lambda F$-cubes must be horizontal and so this must
be zero. Consider next a Haiman arrow $(\rho ,1,0)\to (\rho ',0,0)$
with $\rho '<\rho$. Then this arrow maybe translated so that it is a
strictly southwest pointing Haiman arrow to an $\lambda F$-cube which
we showed in the previous subsection must be zero\footnote{Are you
still with us dear reader? We are deep in the weeds now, but are
almost done.}.
\end{proof}

By the lemma, we conclude that the only sections of
$H^{0}(B,N_{B/S}(\rho x))$ or $H^{0}(B,N_{B/T}(\rho x))$ which
correspond to non-zero Haiman arrows vanish to order at least $\rho$
at $x$, and thus they are necessarily in the image of the maps
\begin{align*}
H^{0}(B,N_{B/S})&\to H^{0}(B,N_{B/S}(\rho x)),\\
H^{0}(B,N_{B/T})&\to H^{0}(B,N_{B/T}(\rho x)).
\end{align*}
By Lemma~\ref{lem: h0 of normal bundles of B in S and T},
$H^{0}(B,N_{B/S})=0$. On the other hand, $H^{0}(B,N_{B/T})$ can be
non-zero and these deformations do occur, they correspond to global
deformations of $B$ in the $K_{S}$-direction.

In conclusion, we have completely classified all possible Haiman
arrows up to equivalence and have thus constructed an explicit basis
for 
\[
\Hom (I_{C},\O_{C})\cong \Ker (\Phi_{1}\oplus \Phi_{2}).
\]
They consist of the $2|\lambda |-\lambda_{1}$ Haiman arrows to
$\lambda F$-cubes which do not go to $B$-cubes and the
$h^{0}(B,N_{B/T})=h^{0}(B,N_{B/X})$ dimensional space of Haiman arrows
going to $B$-cubes. We have thus proved that
\[
\dim \Hom (I_{C},\O_{C}) = h^{0}(B,N_{B/T}) + 2|\lambda |-\lambda_{1}
\]
for $C=B\cup \lambda F$. Our argument extends essentially word for
word to the case where $C=B\cup_{i}(\lambda^{(i)}F_{x_{i}})$ has
several partition thickened fibers. Whether the fiber is smooth or
nodal plays no role. We have thus proved Theorem~\ref{thm: Ext
computation}. \qed 

\subsection{Proof of Theorem~\ref{thm: locus of comb curves is smooth}}\SubSecSpace 

Let $C=B\cup_{i=1}^{n}\left(\lambda^{(i)}F_{x_{i}} \right)$ be a
partition thickened comb curve and let 
\[
d=\sum_{i=1}^{n} |\lambda^{(i)}|,\quad l=\sum_{i=1}^{n}\lambda_{1}^{(i)}.
\]
To prove Theorem~\ref{thm: locus of comb curves is smooth} it will
suffice to construct a flat family of distinct subschemes of $X$,
containing $C$ as a member, and over a base $W$ which is smooth and of
dimension
\[
h^{0}(N_{B/T}) + 2d - l.
\]
Indeed, Theorem~\ref{thm: Ext computation} then implies that the
induced injective map $W\to \Hilb (X)$ is a local isomorphism and the
assertion of Theorem~\ref{thm: locus of comb curves is smooth}
follows.

Let 
\[
H^{0}:=H^{0}(B,N_{B/T}),
\]
and let 
\[
V_{l}\subset \Hilb^{d}(T)
\]
be the stratum defined in Theorem~\ref{thm: strata of Hilb(C^2) with
fixed intersection is smooth}. Let $W=H^{0}\times V_{l}$ so by
Theorem~\ref{thm: strata of Hilb(C^2) with fixed intersection is
smooth}, $W$ is smooth and of dimension $h^{0}(N_{B/T})+2d-l$. We wish
to construct a family over $W$ of distinct subschemes.

Since $T=\Tot (N_{B/T})$, given $\theta \in H^{0}$, we get an
automorphism of $T$, which we call $\Theta$, given by
\[
\Theta :(p,v)\mapsto (p,v+\theta (p)),
\]
where $p\in B$ and $v\in
T|_{p}$. 

We will construct a family of subschemes of $X$, flat over the base
$H^{0}\times V_{l}$, which over a point $(\theta ,Z)\in H^{0}\times
V_{l}$ is the subscheme
\[
C_{\theta} = \Theta (B)\cup p^{-1}(\Theta (Z)).
\]
Clearly, all such subschemes are distinct, and moreover, every
partition thickened comb curve is of the above form (with $\theta =0$
so that $\Theta = \id$).

Formally, we construct the universal subscheme 
\[
\mathcal{C}\subset H^{0}\times V_{l}\times X
\]
flat over $H^{0}\times V_{l}$ as follows. Consider the diagram
\[
\begin{diagram}[h=0.75cm]
&&&&&&\\
&&H^{0}\times V_{l}\times X &&&&\\
&& \uTo^{i} \dTo_{p} &&&&\\
&&H^{0}\times V_{l}\times T & \rTo^{\Theta} & H^{0}\times V_{l}\times  T &&\\
&&\dTo_{\pi_{1}} && \dTo_{\pi_{2}} &&\\
\mathcal{B} &\Into & H^{0}\times T && V_{l}\times T& \lInto & \mathcal{Z}.
\end{diagram}
\]
In the above diagram, $\mathcal{Z}\subset V_{l}\times T$ is the family of 
subschemes of $T$ induced by the universal subscheme over
$\Hilb^{d}(T)$, $\mathcal{B}\subset H^{0}\times T$ is the family of
curves in $T$ given by $H^{0}$, explicitly $\mathcal{B}$ is given by
the set of points $(\theta ,p,\theta (p))$. The maps $\pi_{1}$ and
$\pi_{2}$ are the obvious projections and the maps $p$ and $i$ are the
projection and the zero section of the elliptic fibration $X\to T$. We
are also adopting the general abuse of notation that we drop factors
of the identity map from the notation, that is if $f:A\to B$ we denote
also by $f$ the map $f\times \id_{C}:A\times C\to B\times C$.

Then the subscheme
\[
\mathcal{C} = i(\pi_{1}^{-1}(\mathcal{B})) \cup \left(p\circ \Theta
\circ \pi_{2} \right)^{-1}(\mathcal{Z}) \subset H^{0}\times V_{l}\times X
\]
is the desired universal subscheme over $H^{0}\times V_{l}$. \qed

\subsection{Proof of Theorem~\ref{thm: strata of Hilb(C^2) with fixed intersection
is smooth}}\SubSecSpace 

The constructible function $\Hilb^{d}(T)\to \ZZ $ given by
\[
Z\mapsto \length (Z\cap B)
\]
is upper semi-continuous and thus $V_{l}$ is locally closed.

There is a dense open set on $V_{l}$ isomorphic to
\[
\Sym^{l}(B)\times \Hilb^{d-l}(T-B),
\]
which is clearly smooth and of dimension $2d-l$. Therefore, to prove
the theorem if suffices to show that 
\[
\dim T_{[Z]} V_{l} = 2d-l,
\]
where $Z\subset T$ is a subscheme which is set theoretically (but not
necessarily scheme theoretically) supported on $B$. Moreover, we can
easily reduce to the case where $Z$ is supported at a single point
$p\in B$. By choosing formal local coordinates $(r,s)$ on $T$ at $p$
such that $B=\{s=0 \}$, we are reduced to considering the case
\[
T=\Spec \CC [r,s],\quad B=\{s=0 \},\text{ and $Z\subset T$ supported
at $0$.}
\]
Finally, since $\left(\CC^{*} \right)^{2}$ acts on $V_{l}$ in this
case, it suffices to compute $\dim T_{[Z]} V_{l}$ at $\left(\CC^{*}
\right)^{2}$-fixed subschemes $Z\subset T$. Recall that the fixed
subschemes are given by $Z_{\lambda}$ (see Section~\ref{sec: reduction
to thickened comb curves}) defined by monomial ideals
$I_{\lambda}\subset \CC [r,s]$ corresponding to partitions $\lambda$
of $d$ which in this case have $\lambda_{1}=l$, because $\length (Z_{\lambda}\cap B)=l$.

Therefore, we need only prove the following lemma:
\begin{lemma}\label{lem: dim of Vl at Zlambda is 2d-l}
Let $\lambda =(\lambda_{1}\geq \dotsb \geq \lambda_{k})$ be a
partition of $d$ with $\lambda_{1}=l$. Let $Z_{\lambda}\subset
\CC^{2}=\Spec \CC [r,s]$ be defined by the monomial ideal
\[
I_{\lambda} = (r^{\lambda_{1}},r^{\lambda_{2}}s,\dotsc ,r^{\lambda_{k}}s^{k-1},s^{k}).
\]
Let $V_{l}\subset \Hilb^{d}(\CC^{2})$ be as in Theorem~\ref{thm:
strata of Hilb(C^2) with fixed intersection is smooth}. Then
\[
\dim T_{[Z_{\lambda}]}V_{l} = 2d-l.
\]
\end{lemma}
\begin{proof}
The tangent space
\[
T_{[Z_{\lambda}]}V_{l} \subset T_{[Z_{\lambda}]} \Hilb^{d}(\CC^{2})
\]
is cut out by the equations obtained by linearizing the condition
\[
\length (Z\cap \{s=0 \})=l
\]
at $Z_{\lambda}$. In
\cite{Haiman1998}, M.~Haiman has given a very explicit basis for
$T_{Z_{\lambda}}\Hilb^{d}(\CC^{2})$ in terms of what we called Haiman
arrows in Section~\ref{subsec: combinatorics of Haiman
arrows}. Namely, consider all Haiman arrows $(\rho ,\sigma )\to (\rho
',\sigma ')$ of one of the following two forms:
\begin{enumerate}
\item  $(\rho ,\sigma )\to (\rho ',\sigma ')$ is a southeast
pointing arrow with $(\rho ,\sigma )$ located at a box just above the
top of a column of $\lambda$ and $(\rho ',\sigma ')$ located at a box
which is the furthest to the right in a row of $\lambda$.
\item  $(\rho ,\sigma )\to (\rho ',\sigma ')$ is a northwest
pointing arrow with $(\rho ,\sigma )$ located at a box just to the
right of a row of $\lambda$ and  $(\rho ',\sigma ')$ located at a box
which is at the top of a column of $\lambda$.
\end{enumerate}
There are $2d$ such arrows, $d$ of each kind. The infinitesimal
deformation corresponding to an arrow $(\rho ,\sigma )\to (\rho
',\sigma ')$ is given by deforming the element $r^{\rho}s^{\sigma}\in
I_{\lambda}$ to
\[
r^{\rho}s^{\sigma}+\epsilon \, r^{\rho '}s^{\sigma '},
\]
where $\epsilon^{2}=0$.

For each $\phi :(\rho ,\sigma )\to (\rho ',\sigma ')$ let
$I_{Z_{\lambda}}(\phi )$ be the corresponding deformed ideal and
consider
\[
\dim \left(\frac{\CC [r,s] }{I_{Z_{\lambda}}(\phi) + (s)} \right).
\]
If $\sigma '>0$ then $I_{Z_{\lambda}}(\phi ) + (s) =
I_{Z_{\lambda}} + (s)$: Haiman arrows $(\rho ,\sigma )\to (\rho
',\sigma ')$ with $\sigma '>0$ preserve the condition $\length (Z\cap
\{s=0 \})=l$ and hence lie in $T_{[Z_{\lambda}]}V_{l}$.

If $\sigma '=0$ there are two possibilities:
\begin{enumerate}
\item $(\rho ,\sigma )$ is just above a column of $\lambda$ and $(\rho
',\sigma ')=(\lambda_{1}-1,0)$, or
\item $\phi$ is of the form $(\lambda_{1},0)\to (\rho ',0)$ for $0\leq
\rho '<\lambda_{1}$. 
\end{enumerate}
In Case 2, we have 
\[
\frac{\CC [r,s]}{I_{Z_{\lambda}}(\phi ) + (s)} \cong \frac{\CC
[r]}{(r^{\lambda_{1}}+\epsilon r^{\rho '})},
\]
which has dimension $\lambda_{1}=l$ for all values of $\epsilon$ since
$r^{\lambda_{1}}+\epsilon r^{\rho '}$ has degree $\lambda_{1}$ for all
values of $\epsilon$.

In Case 1, we have
\[
\frac{\CC [r,s]}{I_{Z_{\lambda}}(\phi ) + (s)} \cong \frac{\CC
[r]}{(r^{\lambda_{1}},\epsilon  r^{\lambda_{1}-1})}
\]
which, for non-zero values of $\epsilon$, has dimension
$\lambda_{1}-1$.

Thus we have found that $T_{[Z_{\lambda}]}V_{l}$ is spanned by all the
arrows in the Haiman basis except for the $\lambda_{1}=l$ arrows given
by Case 1 above and therefore $\dim T_{[Z_{\lambda}]}V_{l} = 2d-l$.
\end{proof}

\presectionspace
\appendix
\section{Odds and Ends}\label{appendix: odds and ends}

\subsection{Elliptic surfaces}\label{subsec: elliptic surfs}\SubSecSpace

Let $p : S \rightarrow B$ be a non-trivial elliptic surface with
section $B \subset S$. For simplicity, we assume that all singular
fibers are irreducible nodal rational curves.

Let $X=\Tot (K_{S})$ and let $T=\Tot (K_{S}|_B)$, then clearly we have
\[
N_{B/X}\cong N_{B/S}\oplus N_{B/T}. 
\]

\begin{lemma}\label{lem: h0 of normal bundles of B in S and T}
$h^{0}(N_{B/S})=0$ and $h^{0}(N_{B/T})=\chi (\O_{S})-\chi (\O_{B})$. 
\end{lemma}
\begin{proof}
By a well known fact about elliptic surfaces (see \cite{Fr-Mo} or
\cite[III.1.1]{Miranda}),
\[
K_{S}\cong p^{*}(K_{B}\otimes L),
\]
where
\[
L^{\vee} = R^{1}p_{*}\O_{S}.
\]
Consequently, $c_{1}(K_{S})^{2}=0$ and so Hirzebruch-Riemann-Roch says
\[
\chi (\O_{S}) = \frac{e(S)}{12}>0,
\]
where positivity of $e(S)$ follows by pushing forward the Euler
characteristic measure on $S$ to $B$
\[
e(S) = \int_{S}\, de =  \int_{B} p_{*}(1)\,de = \text{\# of singular fibers.}
\]
On the other hand
\begin{align*}
\chi (\O_{S})& = \chi (R^{\bullet}p_{*}\O_{S})\\
&=\chi (\O_{B}) - \chi (L^{\vee})\\
&=\deg (L). 
\end{align*}
By adjunction
\[
N_{B/S} \cong (K_{S}^{\vee}|_{B})\otimes K_{B} \cong L^{\vee}.
\]
Thus $\deg (N_{B/S})=\deg (L^{\vee})=-\chi (\O_{S})<0$ and so
$h^{0}(N_{B/S})=0$.

Since 
\[
N_{B/T} = K_{S}|_{B} = K_{B}\otimes L,
\]
we see that
\[
h^{1}(N_{B/T}) = h^{1}(K_{B}\otimes L) = h^{0}(L^{\vee}) =
h^{0}(N_{B/S}) = 0,
\]
and therefore
\begin{align*}
h^{0}(N_{B/T})& = \chi (N_{B/T})\\
&=\deg (K_{B}) +\deg (L)+1-g(B)\\
&=\chi (\O_{S})+g(B)-1\\
&=\chi (\O_{S}) - \chi (\O_{B}). \qedhere
\end{align*}
\end{proof}

By our assumption that $S$ is not a product, 
$$
p^* : \Pic^0(B) \stackrel{\cong}{\longrightarrow} \Pic^0(S)
$$
is an isomorphism \cite[VII.1.1]{Miranda}. For any $\beta \in H_2(S)$,
we denote by $\Hilb^\beta(S)$ the Hilbert scheme of effective divisors
on $S$ in class $\beta$.

Denote by $B \in H_2(S)$ the class of the section $B \subset S$ and by
$F \in H_2(S)$ the class of the fiber. Then we have the following
commutative diagram

\[
\begin{diagram}[h=0.9cm]
\Sym^d(B)       &\rTo     & \Pic^{d}(B)\\
\dTo_{p^{*}}    &         &\dTo_{p^{*}}^{\cong }\\
\Hilb^{dF}(S)   &\rTo     &\Pic^{dF}(S) \\
\dTo_{+B}       &         &\dTo_{\otimes \O_S(B)}^{\cong}\\
\Hilb^{B+dF}(S) &\rTo     &\Pic^{B+dF}(S). 
\end{diagram}
\]

The horizontal arrows are Abel-Jacobi maps. The vertical arrows are
induced by pull-back and adding the section $B \subset S$.

\begin{lemma} \label{lem: Sym(B) = Hilb(S)}
The above maps induce a bijective morphism
$$
\Sym^d(B) \to  \Hilb^{B+dF}(S).
$$
\end{lemma}

\begin{proof}
Clearly $p^*$ gives an isomorphism $\Sym^d(B) \cong \Hilb^{dF}(S)$ and
$+B$ defines a closed embedding $\Hilb^{dF}(S) \hookrightarrow
\Hilb^{B+dF}(S)$. Thus it suffices
to show
$$
\Sym^d(B) \rightarrow \Hilb^{B+dF}(S)
$$ 
is surjective on closed points.

For surjectivity, suppose $D'$ is an effective divisor with class
$B+dF$ which does \emph{not} lie in the image. Firstly, we note that
for any fiber $F$ we have $D' \cdot F = 1$. Therefore $D'$ contains a
section $B' \subset S$ as an effective summand. Moreover $B \neq B'$
or else $D'$ would lie in the image. Next, we take any $D$ in the
image and compare $D$ and $D'$. Then
$$
\O_S(D-D') \in \Pic^0(S) \cong \Pic^0(B).
$$ 
Therefore after re-arranging we find that there are distinct fibers $F_{x_i}$, $F_{y_j}$ and $a_i \geq 0$, $b_j \geq 0$ such that 
$$
B + \sum_i a_i F_{x_i} \sim_{\mathrm{lin}} B' + \sum_j b_j F_{y_j},
$$
where $\sim_{\mathrm{lin}}$ denotes linear equivalence. Hence there exists a pencil $\{C_t \}_{t \in \PP^1}$ of effective divisors such that
$$
C_0 = B + \sum_i a_i F_{x_i}, \ C_{\infty} = B' + \sum_j b_j F_{y_j}.
$$
Now fix a smooth fiber $F$. Then $C_t \cdot F = 1$ for any $t \in \PP^1$, so we get a morphism
$$
\PP^1 \longrightarrow F, \ t \mapsto C_t \cap F.
$$
But $F$ is a smooth elliptic curve so this map is constant. We conclude
$$
B \cap F = C_0 \cap F = C_{\infty} \cap F = B' \cap F.
$$
Since $F$ was chosen arbitrary, we deduce that $B = B'$ which is a contradiction.
\end{proof}

\subsection{Weighted Euler characteristics of symmetric products} \label{power}\SubSecSpace 

In this section we prove the following formula for the weighted Euler
characteristic of symmetric products.

\begin{lemma}\label{lem: formula for euler char of sym products}
Let $B$ be a scheme of finite type over $\CC $ and let $e (B)$ be its
topological Euler characteristic. Let $g:\ZZ _{\geq 0}\to \ZZ
(\!(p)\!)$ be any function with $g (0)=1$. Let $G:\Sym ^{d} (B)\to \ZZ
(\!(p)\!)$ be the constructible function defined by
\[
G (\bolda \boldx )=\prod _{i}g (a_{i}),
\]
for all $\bolda \boldx  = \sum_{i}
a_{i}x_{i} \in \Sym^d(B)$ where $x_i \in B$ are distinct closed points. Then
\[
\sum _{d=0}^{\infty } q^{d} \int _{\Sym ^{d} (B)} G \, de =
\left(\sum _{a=0}^{\infty }g (a) q^{a} \right)^{e (B)}.
\]
\end{lemma}

\begin{remark} \label{MacD}
In the special case where $g=G\equiv  1$, the lemma recovers
MacDonald's formula $$\sum _{d=0}^{\infty }e (\Sym ^{d} (B)) \, q^{d} =
\frac{1}{(1-q)^{e (B)}}.$$ 

The lemma is essentially a consequence of the existence of a power
structure on the Grothendieck group of varieties defined by
symmetric products and the compatibility of the Euler characteristic
homomorphism with that power structure. For convenience, we provide a
direct proof here.
\end{remark}
\begin{proof}
The $d$th symmetric product admits a stratification with strata
labelled by partitions of $d$. Associated to any partition of $d$ is a
unique tuple $(m_{1},m_{2},\dots )$ of non-negative integers with
$\sum _{j=1}^{\infty }j m_{j}=d$. The stratum labelled by
$(m_{1},m_{2},\dots )$ parameterizes collections of points where there
are $m_{j}$ points of multiplicity $j$. The full stratification is
given by
\[
\Sym ^{d} (B) = \bigsqcup_{\begin{smallmatrix} (m_{1},m_{2},\dots )\\
\sum _{j=1}^{\infty }j m_{j}=d  \end{smallmatrix}} \left\{\left(\prod _{j=1}^{\infty }B^{m_{j}} \right) -\Delta  \right\}/ \prod _{j=1}^{\infty }\sigma _{m_{j}}, 
\]
where by convention, $B^{0}$ is a point, $\Delta $ is the large
diagonal, and $\sigma _{m}$ is the $m$th symmetric group. Note that
the function $G$ is constant on each stratum and has value $\prod
_{j=1}^{\infty }g (j)^{m_{j}}$. Note also that the action of $\prod
_{j=1}^{\infty }\sigma _{m_{j}}$ on each stratum is free. 

For schemes over $\CC $, topological Euler characteristic is additive
under stratification and multiplicative under maps which are
(topological) fibrations. Thus
\[
\int _{\Sym ^{d} (B)} G\,\, de = \sum _{\begin{smallmatrix}(m_{1},m_{2},\dots )\\
\sum _{j=1}^{\infty }j m_{j}=d   \end{smallmatrix}} \left(\prod _{j=1}^{\infty } g (j)^{m_{j}} \right) \frac{e (B^{\sum _{j}m_{j}}-\Delta )}{m_{1}!\, m_{2}!\, m_{3}!\dots }.
\]

For any natural number $N$, the projection $B^{N}-\Delta \to
B^{N-1}-\Delta $ has fibers of the form $B-\{N-1\text{ points}
\}$. The fibers have constant Euler characteristic given by $e (B)-
(N-1)$ and consequently, $e (B^{N}-\Delta )= (e (B)- (N-1))\cdot e
(B^{N-1}-\Delta )$. Thus by induction, we find $e (B^{N}-\Delta ) = e
(B)\cdot (e (B)-1)\cdots (e (B)- (N-1))$ and so 
\[
\frac{e (B^{\sum _{j}m_{j}}-\Delta )}{m_{1}!\,m_{2}!\,m_{3}!\cdots } = \binom{e (B)}{m_{1},m_{2},m_{3},\cdots },
\]
where the right hand side is the generalized multinomial coefficient.

Putting it together and applying the generalized multinomial theorem,
we find
\begin{align*}
\sum _{d=0}^{\infty } q^d \int _{\Sym ^{d} (B)} G\,\,de & = \sum _{(m_{1},m_{2},\dots )} \prod _{j=1}^{\infty } \left(g (j) q^{j} \right)^{m_{j}} \binom{e (B)}{m_{1},m_{2},m_{3},\dots }\\
&=\left(1+\sum _{j=1}^{\infty }g (j) q^{j} \right)^{e (B)},
\end{align*}
which proves the lemma.   
\end{proof}

\bibliography{/Users/jbryan/jbryan/resources/mainbiblio}
\bibliographystyle{plain}

\end{document}

%% file: 3Dpartition.tex
\newcounter{x}
\newcounter{y}
\newcounter{z}

\newcommand\xaxis{210}
\newcommand\yaxis{-30}
\newcommand\zaxis{90}

\newcommand\topside[3]{
  \fill[fill=yellow, draw=black,shift={(\xaxis:#1)},shift={(\yaxis:#2)},
  shift={(\zaxis:#3)}] (0,0) -- (30:1) -- (0,1) --(150:1)--(0,0);
}

\newcommand\leftside[3]{
  \fill[fill=red, draw=black,shift={(\xaxis:#1)},shift={(\yaxis:#2)},
  shift={(\zaxis:#3)}] (0,0) -- (0,-1) -- (210:1) --(150:1)--(0,0);
}

\newcommand\rightside[3]{
  \fill[fill=blue, draw=black,shift={(\xaxis:#1)},shift={(\yaxis:#2)},
  shift={(\zaxis:#3)}] (0,0) -- (30:1) -- (-30:1) --(0,-1)--(0,0);
}

\newcommand\cube[3]{
  \topside{#1}{#2}{#3} \leftside{#1}{#2}{#3} \rightside{#1}{#2}{#3}
}

\newcommand\planepartition[1]{
 \setcounter{x}{-1}
  \foreach \a in {#1} {
    \addtocounter{x}{1}
    \setcounter{y}{-1}
    \foreach \b in \a {
      \addtocounter{y}{1}
      \setcounter{z}{-1}
      \foreach \c in {1,...,\b} {
        \addtocounter{z}{1}
        \cube{\value{x}}{\value{y}}{\value{z}}
      }
    }
  }
}


\begin{tikzpicture}[scale=0.35]

\newcommand{\redbox}[3]{ 
\draw[fill=red, draw=black,shift={(\xaxis:#1)},shift={(\yaxis:#2)},shift={(\zaxis:#3)}]
(0,0) rectangle (1,1);
}
\newcommand{\yellowbox}[3]{ 
\draw[fill=yellow, draw=black,shift={(\xaxis:#1)},shift={(\yaxis:#2)},shift={(\zaxis:#3)}]
(0,0) rectangle (1,1);
}
\newcommand{\bluebox}[3]{ 
\draw[fill=blue, draw=black,shift={(\xaxis:#1)},shift={(\yaxis:#2)},shift={(\zaxis:#3)}]
(0,0) rectangle (1,1);
}

\planepartition{{11,11},{11,11},{11},{1},{1},{1},{1},{1},{1},{1},{1}}

\draw [thick,->] (0,11)-- (0,14);
\draw [thick,->] (1.74,-1.74*0.57735)-- (11,-11*0.57735);
\draw [thick,->]  (-9.5,-9.5*0.57735)-- (-11,-11*0.57735);

\node [right] at (0,13) {$t$};
\node [right] at (11,-11*0.57735) {$s$};
\node [left] at  (-11,-11*0.57735) {$r$};
\node [above] at  (-9,-7*0.57735) {$B$};
\node [right] at (2,10) {$\lambda F$};
\node  at (6,6*0.57735) {$R$};
\node  at (-6,6*0.57735) {$S$};
\node  at (0,-6) {$T$};

\shade [ball color=gray] (-2.6,-.6) circle [radius=0.4cm];
\shade [ball color=gray] (1.8,-1.2) circle [radius=0.4cm];
\shade [ball color=gray] (-0.85,-1.7) circle [radius=0.4cm];

\end{tikzpicture}
